\documentclass[12pt]{article}

\usepackage{amsfonts}
\usepackage{amsmath}
\usepackage{amsthm} %removed this on 12th April because of latex error.
\usepackage{epsfig}
\usepackage{graphics}

\usepackage{latexsym}
\usepackage{amssymb,exscale}
\usepackage[all]{xy}
\usepackage{hyperref}
\usepackage{verbatim} % for commenting out big chunks of text

\usepackage{pslatex}
\theoremstyle{plain} %documentation says there are only three styles
    \newtheorem{theorem}{Theorem}
    \newtheorem{lemma}[theorem]{Lemma}

\theoremstyle{definition} % For roman text in the body
    \newtheorem{definition}[theorem]{Definition}
    
    \newtheorem{result}[theorem]{Result}
    \newtheorem{remark}[theorem]{Remark}
    \newtheorem{example}[theorem]{Example}

\def\Gam{\Gamma}
\def\Del{\Delta}

\def\Ome{\Omega}
\def\alp{\alpha}
\def\bet{\beta}
\def\gam{\gamma}
\def\del{\delta}
\def\eps{\epsilon}

\def\lam{\lambda}

\def\Ome{\Omega}
\def\sig{\sigma}

\def \K{{\mathcal K}}

\def\C{\mathbb{C}}
\def\D{\mathbb{D}}

\def\GG{{\mathcal G}}
\def\H{{\bf{H}}}

\def\M{\Omega}

\def\PP{{\mathcal P}}

\def\R{\mathbb{R}}
\def\Sym{{\mathcal S}}
\def\S{\mathbb{S}}
\def\U{{\mathcal U}}
\def\X{{\mathcal X}}

\def\f{{\bf f}}

\def\q{{\bf q}}

\def\summ{\sum\limits}
\def\intt{\int\limits}
\def\prodd{\prod\limits}

\def\tends{\rightarrow}

\def\sgn{{\mb{sgn}}}
\def\Mid{\left\vert \right.}
\def\l{\left}
\def\r{\right}
\def\<{\langle}
\def\>{\rangle}
\def\hsp{\hspace}
\def\mb{\mbox}

\newcommand{\E}{\mbox{\bf E}}
\def\I#1{{\bf 1}_{#1}}

\def\bar{\overline}
\def\P{{\bf P}}

\def\convd{\stackrel{d}{\rightarrow}}

\def\eqd{\stackrel{d}{=}}
\def\given{\left.\vphantom{\hbox{\Large (}}\right|}

\newcommand\Tr{{\mbox{Tr}}}
\newcommand\mnote[1]{} %off
\newcommand\be{\begin{equation*}}

\newcommand\ee{\end{equation*}}
\newcommand\lb[1]{\label{#1}\mnote{#1}}
\newcommand\ben{\begin{equation}}
\newcommand\een{\end{equation}}
\newcommand\bes{\begin{eqnarray*}}
\newcommand\ees{\end{eqnarray*}}

\newcommand{\sm}{{\raise0.3ex\hbox{$\scriptstyle \setminus$}}}

\def\mb{\mbox}
\def\l{\left}
\def\r{\right}
\def\sig{\sigma}
\def\lam{\lambda}
\def\alp{\alpha}
\def\eps{\epsilon}
\def\tends{\rightarrow}

\renewcommand{\phi}{\varphi}

\def\Kdet{{\mathbb K}}

\author{Manjunath Krishnapur}
\title{From random matrices to \\ random analytic functions}
\begin{document}
 \bibliographystyle{abbrv}
\maketitle
\section{Leading up to the results}
Singular points of random matrix-valued analytic functions are a common generalization of eigenvalues of random matrices and zeros of random polynomials. The setting is that we have an analytic function of $z$ taking values in the space of $n\times n$ matrices. Singular points are those (random) $z$ where the matrix becomes singular, that is, the zeros of the determinant. This notion was introduced in the Ph.D thesis~\cite{krithesis} of the author, where some basic facts were found. Of course, singular points are just the zeros of the (random analytic function) determinant, so in what sense is this concept novel? 
 
In case of random matrices as well as random analytic functions, the following features may be observed.
\begin{enumerate}
\item For very special models, usually with independent Gaussian coefficients or entries, one may solve {\em exactly} for the distribution of zeros or eigenvalues. 
\item For more general models with independent coefficients or entries, under rather weak assumptions on moments, one can usually analyze the empirical measure of eigenvalues or zeros as the size of the matrix increases or the degree of the polynomial goes to infinity. 
\item Substituting independence with assumptions of particular kinds of symmetry and dependence of entries or coefficients, eigenvalues and zeros have been studied with varying degrees of success.
\end{enumerate}
The point here is that the determinant of a random matrix-valued analytic function has coefficients that are dependent in a very complicated way and one might not expect it to be tractable. Nevertheless, 
\begin{itemize}
\item  In this paper we demonstrate that certain special models of random matrix-valued analytic functions based on independent complex Gaussians have singular sets that turn out to be determinantal point processes in the two dimensional sphere or the hyperbolic plane! This is the exactly solvable situation.
\item In a subsequent paper~\cite{kricircularlawgeneralizations}, we shall study asymptotics of the counting measure on the singular set of random matrix-valued analytic functions, under weak assumptions of independence (and some moment conditions),  as the matrix size goes to infinity. This will be a generalization of the circular law for non-Hermitian random matrices with independent entries.
\item Numerous questions suggest themselves, taking unitary or Hermitian matrix coefficients in a random polynomial, for example. Answers are yet to suggest themselves.
\end{itemize}
Despite its natural appeal, the concept of a random matrix-valued analytic function does not seem to have been considered in the literature, perhaps because the focus in random matrix theory
has been mostly on eigenvalues in one dimension (real line or
the circle). We do not know a way to force singular points to lie on the line  (except the case of eigenvalues of Hermitian matrices).

A necessary notion needed to even state some of our results is that of a {\bf determinantal point process}, first defined by Macchi~\cite{macchi}. For the reader not familiar with them, a brief introduction to determinantal processes is given in the Appendix. This is sufficient for the purposes of this paper, but to know  more, the reader may consult the surveys~\cite{soshnikov00} or \cite{hkpv}. The reader interested merely in our results and proofs may now jump directly to the next section. The rest of this section is devoted to motivating the results and establishing the context and is not logically necessary to read the rest of the paper. First some notations.

\vspace{2mm}
\noindent{\bf Notations:} $\D $ is the unit disk in the complex plane. The two-dimensional sphere $\S^2$ will always be identified with $\C\cup \{\infty\}$ via stereographic projection. $m$ denotes Lebesgue measure. $g$ and $g_i$, $i\ge 0$, are always  independent standard complex Gaussian random variables, with density $\pi^{-1}\exp\{-|z|^2\}$ in the plane. $G$ and $G_i$, $i\ge 0$, are $n\times n$ matrices whose entries are i.i.d. standard complex Gaussians. The group of permutations of a set $T$ will be denoted by $\Sym(T)$. When  $T=\{1,2,\ldots ,k\}$ we just write $\Sym_k$. We denote the set $\{1,2,\ldots ,k\}$ by $[k]$. And $\U(N)$ is the group of $N\times N$ unitary matrices.

\vspace{3mm}
The results of this paper were motivated by two well known results, one from the realm of random matrices and another concerning zeros of random analytic functions.  The first is due to Ginibre~\cite{gin}, who found the exact distribution of eigenvalues of a random matrix with independent standard complex Gaussian entries and the second is due to Peres and Vir\'{a}g~\cite{pervir} who discovered the exact distribution of zeros of the random power series with independent standard complex Gaussian coefficients.

\begin{result}[Ginibre(1965)]\label{thm:ginibre}
The set of eigenvalues of the  $n\times n$ matrix $G$ with i.i.d. standard complex Gaussian entries is a determinantal point process with kernel
\begin{equation}\label{eq:gindetform}
  \Kdet(z,w) = \summ_{k=0}^{n-1}
  \frac{(z\bar{w})^k}{k!},
\end{equation}
with respect to the reference measure $d\mu(z)=\frac{1}{\pi}e^{-|z|^2}$. 
The corresponding Hilbert space $\H=\mb{span}\{1,z,\ldots ,z^{n-1}
\}\subset L^2(\C,\frac{e^{-|z|^2}}{\pi}dm(z))$.
\end{result}

\begin{result}[Peres and Vir\'ag(2003)]\label{thm:pervir}
 Let $\f(z)=g_0+g_1z+g_2z^2+\ldots$  be the random analytic function (the radius of convergence is $1$) whose coefficients $g_i$ are i.i.d. standard complex Gaussians. Then the zeros of $\f$ form a determinantal point process on the unit disk $\D$ with the kernel (the Bergman kernel of the unit disk)
 \begin{equation*}
   \Kdet(z,w) = \frac{1}{\pi(1-z\bar{w})^2},
 \end{equation*}
with respect to  the background measure $d\mu(z)=\frac{1}{\pi}dm(z)$ on $\D$.
 The corresponding Hilbert space $\H=\mb{span}\{1,z, z^2, \ldots  
\}\subset L^2(\D,\frac{dm(z)}{\pi})$ is the space of all analytic
functions in $L^2(\D)$.
\end{result}
Knowing that a point process is determinantal greatly facilitates studying its properties. This motivates us to ask whether these two are isolated results or whether they are part of a bigger picture. As a start, let us draw a list of random analytic functions whose zero sets are known to be determinantal.

\begin{itemize}
\item $z-g=0$ : One zero, $g$, with standard complex Gaussian
distribution on $\C$. $\H=\mb{span}\{1\}$ in
$L^2(\C,\frac{e^{-|z|^2}}{\pi}dm(z))$. Determinantal, of course!

\item $zg_1-g_2=0$ : One zero, $\frac{g_2}{g_1}$, that has density $\frac{1}{\pi(1+|z|^2)^2}dm(z)$. This is just the push-forward of the uniform measure on $\S^2=\C \cup\{\infty\}$ under stereographic projection. Again, a one-point process is determinantal. In this case,  $\H=\mb{span}\{1\}$ in $L^2(\C \cup\{\infty\},\frac{1}{\pi(1+|z|^2)^2}dm(z))$. 

\item $g_0+zg_1+z^2g_2+\ldots=0$. This is the i.i.d. power series of Theorem~\ref{thm:pervir}.
\end{itemize}
Our key observation is that Ginibre's result (Result~\ref{thm:ginibre}) describes the law of zeros of the random analytic function $\det(zI-G)$, which may in turn be thought of as the matrix version of the analytic function $z-g$, the first of the three examples above. This suggests that we consider the matrix versions of the second and third examples. This leads us to two families of random matrix-valued analytic functions. 
\begin{enumerate}
\item $zG_1-G_2$, where $G_1,G_2$ are $n\times n$ independent matrices with
  i.i.d. standard complex Gaussian entries.

\item $G_0+zG_1+z^2G_2+\ldots$, where $G_k$ are
  independent $n\times n$ matrices with each $G_k$ having i.i.d. standard complex Gaussian entries.
\end{enumerate}
The analogy with Ginibre's result strongly suggests that the singular points of  these matrix-valued analytic functions might be determinantal point processes.  But how to guess which determinantal processes?

\noindent{\bf Invariance to the rescue:} The key feature that allows us to guess which determinantal processes, is invariance under a large group of transformations. A point process $\X$ on the space $\M$ is said to be invariant (in distribution) under a transformation $T:\M\tends \M$ if $T(\X)\eqd \X$.

\begin{itemize}
\item  First consider the matrix-valued analytic function $zG_1-G_2$. We claim that the singular points (which are just the eigenvalues of $G_1^{-1}G_2$) are invariant under all linear fractional transformations
\be
 \lam \tends \frac{\lam\alp- \bar{\bet}}{\lam \bet-\bar{\alp}}
\ee
with $\alp,\bet$ complex numbers satisfying $|\alp|^2+|\bet|^2=1$. These are precisely rotations of the two-dimensional sphere, when the sphere is identified with $\C\cup \{\infty\}$ via stereographic projection.

To see this, let $\alp,\bet$ be complex numbers such that $|\alp|^2+|\bet|^2=1$. Then define 
\be
 C=\alp G_1-\bet G_2 \mb{ and } D=\bar{\bet} G_1 +\bar{\alp} G_2.
\ee
Then $(C,D)$ has the same distribution as $(G_1,G_2)$. This implies that the solutions to $\det(zC-D)=0$ have the same distribution as the solutions to $\det(zG_1-G_2)$. On the other hand,
\bes
 \det(zC-D) &=& \det\l((z\alp- \bar{\bet})G_1-(z\bet-\bar{\alp})G_2 \r) \\
            &=& (z\bet-\bar{\alp})^n \det\l(\frac{z\alp- \bar{\bet}}{z\bet-\bar{\alp}}G_1-G_2 \r).
\ees
Thus the zeros of $\det(zC-D)$ are precisely $\l\{\frac{\lam_i\alp- \bar{\bet}}{\lam_i\bet-\bar{\alp}}\r\}$ where $\{\lam_i\}$ are the zeros of $\det(zG_1-G_2)$. Thus we have
\ben \lb{eq:sphericalinvariance}
 \l\{\frac{\lam_i\alp- \bar{\bet}}{\lam_i\bet-\bar{\alp}}\r\}_{1\le i\le n} \eqd \{\lam_i\}_{1\le i \le n}
\een
for any $\alp,\bet$ with $|\alp|^2+|\bet|^2=1$ which is what we claimed.

\item Next  consider the matrix-valued analytic function $M(z):=\summ_{k=0}^{\infty} G_k z^k$ where $G_k$ are independent random matrices with i.i.d. standard complex Gaussian entries. We claim that the set of singular points is invariant in distribution under the isometries of the hyperbolic plane, namely the linear fractional transformations
\be
 \phi(\lam) = \frac{\alp \lam + \bet}{\bar{\bet}\lam+\bar{\alp}}, \hsp{2mm} |\alp|^2-|\bet|^2=1,
\ee
that map the unit disk injectively onto itself. These are precisely the conformal automorphisms of the unit disk. To see this, observe that for each $1\le i,j \le n$, 
\be
 M_{i,j}(z) = G_0(i,j)+G_1(i,j)z+G_2(i,j)z^2+\ldots 
\ee
is a copy i.i.d. power series of Theorem~\ref{thm:pervir}, and  the $M_{i,j}(\cdot)$s are themselves independent random functions. It is known (see \cite{ST1}) that
\ben \label{eq:invarianceforiidps}
 M_{i,j}(\phi(\cdot)) \eqd \phi'(\cdot)^{-\frac{1}{2}}M_{i,j}(\cdot).
\een
This is because $M_{1,1}$ is the Gaussian element of the Hilbert space of analytic functions with the boundary inner product $\int_0^{2\pi} f(e^{i\theta}) \bar{g}(e^{i\theta}) \frac{d\theta}{2\pi}$ (We define this inner product for analytic functions that extend continuously to the boundary, and take the completion). On this Hilbert space,  $f \tends \phi'(\cdot)^{-\frac{1}{2}}f(\phi(\cdot))$ is a unitary transformation. When a unitary transformation is applied to the Gaussian element of the Hilbert space, we get again the  Gaussian element of the Hilbert space, yielding  (\ref{eq:invarianceforiidps})  (of course, the Gaussian element is itself not an element of the Hilbert space, almost surely).

An alternate way is to just check that the centered Gaussian processes on the two sides of (\ref{eq:invarianceforiidps}) have the same covariance kernel $\frac{1}{1-\phi(z)\bar{\phi(w)}}$. The reader may refer to the paper of Sodin and Tsirelson~\cite{ST1} for a more detailed proof. Anyhow, the independence of distinct $M_{i,j}$ and the fact that $\phi'$ is non-random shows that
\be
 \det\l(M(\phi(\cdot))\r) \eqd \phi'(\cdot)^{-\frac{n}{2}} \det(M(\cdot)).
\ee
Now, $\phi'$ is a nowhere vanishing analytic function on the unit disk. Therefore the above equation shows that the singular set of $M(\cdot)$ is invariant in distribution under the action of hyperbolic isometries.
\end{itemize}
These are two special cases of a large class of invariant zero sets introduced in \cite{krithesis}. The general situation is that one applies a homogeneous polynomial of several complex variables (``$\det$'' in our case) to a bunch of i.i.d. copies of a Gaussian analytic function ($zg_1-g_2$ or $g_0+zg_1+z^2g_2+\ldots$ in the two cases). If the individual Gaussian analytic functions have invariant zero sets, then so will the homogeneous polynomial of copies of them. The idea is that constructions (applying a homogeneous polynomial to i.i.d. copies of an analytic function) which are simple in terms of functions, are not simple at all at the level of zeros, and may give something drastically new. 

Let us return to our original question which led to a digression into the issue of invariance. This was the question of guessing what determinantal processes might the singular sets of the two families of random matrix-analytic functions ($zG_1-G_2$ and $G_0+zG_1+z^2G_2+\ldots $) be. We now have obtained a strong restriction: {\em the determinantal process better be isometry-invariant in $\S^2$ or $\D$, respectively}. Further, from the fact that we are looking at zeros of analytic functions, we expect that these determinantal processes are defined by Hilbert spaces of analytic functions. Such determinantal processes were classified in \cite{krithesis} (see Theorem 3.0.5 therein). From this classification, we get the following invariant determinantal processes as the only possible candidates. These processes were, in fact, studied first by Caillol~\cite{caillol} under the name ``one component plasma on the sphere'' and  by Jancovici and T{\'e}llez~\cite{jantell} on the hyperbolic plane.
\begin{enumerate}
\item On the sphere ($\C\cup \{\infty\}$) we have for each $n\in \{1,2,\ldots \}$ an invariant determinantal point process with kernel
\be
 \Kdet(z,w)  = (1+z\bar{w})^{n-1}
\ee
with respect to the background measure $d\mu_n(z) = \frac{n}{\pi(1+|z|^2)^{n+1}}dm(z)$. Invariance of the point process under an analytic transformation $T$  is equivalent to saying that the joint intensities (correlation functions) with respect to Lebesgue measure satisfy
\be
\rho_k(Tz_1,\ldots ,Tz_k) |T'(z_1)|^2\ldots |T'(z_k)|^2 = \rho_k(z_1,\ldots ,z_k).
\ee
For the case at hand, this is easily checked from the fact
\be
\frac{T'(z) \bar{T'(w)}}{(1+T(z)\bar{T(w)})^2}  = \frac{1}{(1+z\bar{w})^2}.
\ee
The parameter $n$ is  the total number of points in the point process, or equivalently, it denotes the first intensity of the point process with respect to the spherical area measure $\frac{1}{\pi(1+|z|^2)^2}$.

\item On the unit disk, we have for each $n>0$, an invariant determinantal point process  with kernel
\be
 \Kdet(z,w)  = \frac{1}{(1-z\bar{w})^{n+1}}
\ee
with respect to the background measure $d\mu_n(z) = \frac{n}{\pi}(1-|z|^2)^{n-1}dm(z)$. Again it is easy to check that these are invariant, now using
\be
\frac{\phi'(z) \bar{\phi'(w)}}{(1-\phi(z)\bar{\phi(w)})^2}  = \frac{1}{(1-z\bar{w})^2}.
\ee
The parameter $n$ is denotes the first intensity of the point process with respect to the hyperbolic measure $\frac{1}{\pi(1-|z|^2)^2}$.
\end{enumerate}
Thus, on each of the sphere and the disk, we have a family of
invariant singular sets and a family of invariant determinantal processes. Then by comparing the first intensities of these determinantal processes and the set of singular points of our matrix-analytic functions, we match the singular sets to determinantal processes.

\section{Statements of results}
We now state our results.
\begin{theorem}\label{thm:sphere} Let $G_1,G_2$ be i.i.d. $n\times n$ matrices with
i.i.d. standard complex Gaussian entries. The zeros of $\det(zG_1-G_2)$ form a
determinantal point process on $\S^2$ with kernel
\begin{equation*}
  \Kdet(z,w) = (1+z\bar{w})^{n-1}
\end{equation*}
with respect to the background measure $d\mu_n(z)=\frac{n}{\pi}\frac{dm(z)}{(1+|z|^2)^{n+1}}$.  
Equivalently, we may say that the defining Hilbert space is the subspace of {\it analytic} functions in
$L^2 (\C \cup \{\infty\}, \mu_n)$.
\end{theorem}

\begin{theorem}\label{thm:disk} Let $G_k$ be i.i.d. $n\times n$ matrices with
i.i.d. standard complex Gaussian entries. Then for each $n\ge 1$, the zeros of $\det\l( G_0+zG_1+z^2G_2+\ldots \r)$ form a determinantal point process on $\D$
with kernel
\begin{equation*}
  \Kdet(z,w) = \frac{1}{(1-z\bar{w})^{n+1}}  
\end{equation*}
with respect to the background measure $d\mu_n(z) = \frac{n}{\pi}
 (1-|z|^2)^{n-1}dm(z)$. Equivalently, we may say that the defining Hilbert space is the subspace of {\it analytic} functions in $L^2 (\D, \mu_n)$.
\end{theorem}
Theorem~\ref{thm:sphere} is proved section~\ref{sec:sphere} via the Schur decomposition of the matrix $G_1^{-1}G_2$, along the lines of Ginibre's proof of Theorem~\ref{thm:ginibre}. Theorem~\ref{thm:disk} will be proved in section~\ref{sec:disk} as a corollary of the following more general theorem which appears to be of potential interest beyond the specific application to Theorem~\ref{thm:disk}. Theorem~\ref{thm:general} is proved in section~\ref{sec:general}.

\begin{theorem} \label{thm:general}
Let $A_N$ be $n\times n$  matrices such that $\sqrt{N}A_N\convd X_0$ for some random matrix $X_0$. Independently of $A$, pick $P,Q$ independent matrices chosen from Haar measure on $\U(N)$ and define the $N\times N$ matrix $V$ by
\ben \label{eq:lawofV}
V = Q^*\l[\begin{array}{cc} 
      A_N & 0 \\
      0 & I_{N-n} \end{array} \r]P^*.
\een
Set $\f_N(z) := \frac{\det(zI+V)}{\det(I+zV^*)}$. Let $X_0$ and $G_i$, $i\ge 1$, be independent random matrices, where $G_i$ have independent standard complex Gaussian entries. Then
\be
 N^{n/2}\f_N(z) \convd \det\l(X_0+\summ_{k\ge 1} G_k z^k\r)
\ee
in the sense that any finite set of coefficients in the power series expansion of $N^{n/2}\f_N(\cdot)$ converge jointly in distribution to the corresponding vector of coefficients in the power series expansion of the right hand side.
%Here the convergence in distribution is with respect to the %Caratheodary topology\footnote{The topology of uniform %convergence on compact sets. Equivalent to convergence of %coefficients in the power series. Implies convergence of zeros.} %on the space of analytic functions on $\D$.
\end{theorem}
The relevance of this theorem to Theorem~\ref{thm:disk} is through a result of {\.Z}yczkowski and Sommers~\cite{zycsom} who  found random matrix models whose eigenvalue distributions are determinantal processes with kernels that are truncated versions of the kernels in Theorem~\ref{thm:disk}. Theorem~\ref{thm:general} gives the distribution of the limiting random analytic function, as the matrix size increases, while the result of {\.Z}yczkowski and Sommers gives the limiting distribution of zeros. Putting the two together we get the distribution of zeros of the limiting random analytic function.

\begin{remark} The statements of Theorem~\ref{thm:sphere} and Theorem~\ref{thm:disk} may already be found in the thesis \cite{krithesis}. Theorem~\ref{thm:sphere} appeared there with a proof but is being published here for the first time. Theorem~\ref{thm:disk} was conjectured in \cite{krithesis} and a partial proof was given, showing that the first and second joint intensities (correlation functions) of the singular set of $G_0+zG_1+z^2G_2+\ldots $ are as claimed. While in \cite{krithesis} we tried to prove Theorem~\ref{thm:disk} by starting with the random matrix-valued analytic function and then finding the distribution of its zeros, in contrast, in this paper we take the opposite direction.
\end{remark}
 
\begin{remark} It is natural to ask whether there are other (perhaps even many) singular sets that are also determinantal. Without claiming that there are not, we would like to emphasize that the determinantal processes in Theorem~\ref{thm:ginibre} together with those in  Theorems~\ref{thm:sphere},\ref{thm:disk}, are {\em the most natural} determinantal point processes in the three canonical surfaces of constant curvature, namely, the plane, the sphere and the hyperbolic plane, respectively. As remarked earlier, these processes were studied first by Caillol~\cite{caillol} on the sphere and  by Jancovici and T{\'e}llez~\cite{jantell} on the hyperbolic plane (related "two-component plasmas" were studied by Forrester, Jancovici and Madore~\cite{forjanmad}).

Independently, in \cite{krithesis}, motivated by an analogous theorem of Sodin~\cite{sodin} for zeros of Gaussian analytic functions, it was proved that on each of these three domains, there is exactly a one parameter family of invariant determinantal point processes  that arise from Hilbert spaces of analytic functions. There are some additional conditions, see Theorem~3.0.5 in \cite{krithesis} for precise statements. These determinantal processes are exactly those that correspond to Hilbert spaces of analytic functions on these domains with respect to the following measures.
\begin{itemize}
\item $d\mu_{\alp}(z)=\frac{\alp}{\pi}e^{-\alp|z|^2} dm(z)$, $\alp>0$ in the plane. In this case, changing $\alp$ merely has the effect of scaling the plane, and therefore these determinantal processes should be thought of as identical. 
\item $d\mu_{\alp}(z) = \frac{\alp}{\pi(1+|z|^2)^{\alp+1}}dm(z)$ for $\alp\in \{1,2,3,\ldots \}$ in the sphere ($\C\cup \{\infty\}$).
\item $d\mu_{\alp}(z) = \frac{\alp}{\pi}(1-|z|^2)^{\alp-1}dm(z)$ for $\alp>0$ in the unit disk.
\end{itemize}
The determinantal processes appearing in Theorems~\ref{thm:sphere},\ref{thm:disk} are precisely these canonical ones, while the determinantal processes in Ginibre's theorem converge (as the matrix size increases) to the canonical determinantal point process on the plane. The upshot of all this is that determinantal singular sets are somewhat special and may not be all that abundant. 

Note that while canonical determinantal processes in the unit disk exist for every $\alp>0$, Theorem~\ref{thm:disk} gives a singular-set interpretation only for positive integer values of $\alp$.
\end{remark}

\section{Spherical ensembles}\label{sec:sphere}
 Let $\X$ denote the set of singular points of $zG_1-G_2$. Since the number of points is exactly $n$, Theorem~\ref{thm:sphere} is equivalent (see the facts stated after definition~\ref{def:det} in the appendix) to saying that the joint density of the singular points is proportional to 
\be
 \Mid \Del(z_1,\ldots ,z_n) \Mid^2 \prodd_{k=1}^n \frac{1}{(1+|z_k|^2)^{n+1}}.
\ee
(If $\{P_1,\ldots ,P_n\}$ are the points on the two-dimensional sphere obtained by stereographic projection of $z_1,\ldots ,z_n$, then the density of these points with respect to Lebesgue measure on $(\S^2)^n$ is simply
\be
 \prodd_{i<j} \|P_i-P_j\|_{\R^3}^2
\ee
where $\|\cdot \|_{\R^3}$ is the Euclidean norm in $\R^3$). 

The following lemma will greatly simplify the job of integrating out auxiliary variables later.
\begin{lemma}\label{lem:invdist} Let $\X$ be a point process on $\C$
  with $n$ points almost surely. Assume that the $n$-point correlation
  function (equivalently the density) of $\X$ has the form
\[ p(z_1,\ldots ,z_n) = \Mid \Del(z_1,\ldots ,z_n) \Mid^2 V(|z_1|^2,\ldots ,|z_n|^2). \]
Here $\Del(z_1,\ldots ,z_n)$ denotes the Vandermonde factor $\prodd_{i<j}(z_j-z_i)$. 

Suppose also that $\X$ has a distribution invariant under automorphisms of the sphere $\S^2$, i.e.,  under the transformations $\phi_{\alp,\bet}(z)=\frac{\alp z +\bet}{-\bar{\bet}z+\bar{\alp}}$, for any $\alp,\bet$ satisfying $|\alp|^2+|\bet|^2=1$.  
Then 
\begin{equation}\label{eq:specialformofdensity}
V(|z_1|^2,\ldots ,|z_n|^2)=\mb{Const.} \prodd_{k=1}^n
\frac{1}{(1+|z_k|^2)^{n+1}}.  
\end{equation} 
\end{lemma}

\begin{proof}[{\bf Proof of Lemma~\ref{lem:invdist}}] The claim is that the probability density of the $n$
  points of $\X$ (in exchangeable random order) with respect to Lebesgue measure is 
  \begin{equation*}
    q(z_1,\ldots,z_n) :=\mb{Const.} \Mid \Del(z_1,\ldots ,z_n) \Mid^2
\prodd_{k=1}^n \frac{1}{(1+|z_k|^2)^{n+1}}.
  \end{equation*}
First let us check that the density $q$ is invariant under the
isometries of $\S^2$. For this let $\phi_{\alp,\bet}(z)=\frac{\alp z
  +\bet}{-\bar{\bet}z+\bar{\alp}}$, with $\alp,\bet$ satisfying
$|\alp|^2+|\bet|^2=1$.  Then, 
\begin{equation}\label{eq:derphi}
\phi^{\prime}(z)= \frac{1}{(-{\bar \bet}z +{\bar \alp})^2}.
\end{equation}
\begin{equation}\label{eq:metphi}
1+\Mid \phi(z) \Mid^2 = \frac{1+\Mid z \Mid^2}{\Mid -{\bar \bet}z +{\bar \alp}\Mid^2}.
\end{equation}
\begin{equation}\label{eq:diffphithis}
\phi(z) - \phi(w) = \frac{z-w}{(-{\bar \bet}z +{\bar \alp})(-{\bar \bet}w +{\bar \alp})}.
\end{equation}
From (\ref{eq:derphi}),(\ref{eq:metphi}) and (\ref{eq:diffphithis}), it
follows that
\begin{equation}\label{eq:invforq}
 q\l(\phi(z_1),\ldots ,\phi(z_n)\r)\prodd_{k=1}^n |\phi^{\prime}(z_k)|^2 =
 q\l(z_1,\ldots, z_n\r), 
\end{equation}
which shows the invariance of $q$.  

Invariance of $\X$ means that $\forall \alp,\bet$ with
$|\alp|^2+|\bet|^2=1$, and for every $z_1,\ldots ,z_n$, we have  (with $\phi=\phi_{\alp,\bet}$)
\begin{equation}\label{eq:invforp}
 p\l(\phi(z_1),\ldots ,\phi(z_n)\r)\prodd_{k=1}^n |\phi^{\prime}(z_k)|^2 =
 p\l(z_1,\ldots, z_n\r). 
\end{equation}
Set $W(z_1,\ldots ,z_n)=\frac{p(z_1,\ldots
  ,z_n)}{q(z_1,\ldots,z_n)}$. Then,  we get
\begin{itemize}
\item $W\l(z_1,\ldots ,z_n\r)$ is a function of $|z_1|^2,\ldots ,|z_n|^2$  only, by the assumption on $p$ and the definition of $q$.
\item $W\l(\phi(z_1),\ldots ,\phi(z_n)\r)=W\l(z_1,\ldots ,z_n \r)$ for
  every $z_1,\ldots ,z_n$ from (\ref{eq:invforp}) and
  (\ref{eq:invforq}). 
\end{itemize}   
We claim that these two statements imply that $W$ is a constant. To
see this fix $z_k=r_k e^{i\theta_k}$, $1\le k\le n$, such that $r_1
<r_k$ for $k\ge 2$. Let $\alp=\frac{1}{\sqrt{1+r_1^2}},
\bet=-\frac{z_1}{\sqrt{1+r_1^2}}$. Then $|\alp|^2+|\bet|^2=1$ and so
$\phi_{\alp,\bet}$ is an isometry of $\S^2$. From the above stated
properties of $W$, we deduce,
\begin{eqnarray*}
  W(z_1,\ldots,z_n) &=& W\l(\phi(z_1),\ldots,\phi(z_n)\r) \\
                    &=& W\l(0,\frac{z_2-z_1}{1+z_2{\bar z_1}},\ldots
                    ,\frac{z_n-z_1}{1+z_n{\bar z_1}}\r) \\
                    &=& W\l(0, \given
                    \frac{r_2e^{i\theta_2}-z_1}{1+r_2e^{i\theta_2}{\bar
                    z_1}} \given, \ldots
                    , \given \frac{r_ne^{i\theta_n}-z_1}{1+r_ne^{i\theta_n}{\bar
                    z_1}} \given \r). 
\end{eqnarray*}
Take $z_1=1$ and $1<r_k<1+\eps$. Then as $\theta_k$, $2\le k\le n$ vary
independently over $[0,2\pi]$, the quantities $\Mid
                    \frac{r_ke^{i\theta_k}-z_1}{1+r_ke^{i\theta_k}{\bar
                    z_1}} \Mid$ vary over the intervals
                $\l[\frac{r_k-1}{r_k+1},\frac{r_k+1}{r_k-1}\r]$. However the left side, $W(z_1,\ldots ,z_n)$ does not change because $W$ is a function of $r_k$s only. By
                our choice of $r_k$s, this means that
                \begin{equation*}
                  W(0,t_2,\ldots, t_n) = \mb{Constant} \hsp{1cm}
                  \forall t_k\in \l[\eps,\frac{1}{\eps}\r].
                \end{equation*}
$\eps$ is arbitrary, hence $W(0,t_2,\ldots,t_n)$ is
constant. This implies that $W(0,z_2,\ldots ,z_n)$ is constant and therefore $W(z_1,\ldots ,z_n)$ is constant. 

This shows that $p(z_1,\ldots, z_n)=\mb{Const.} q(z_1,\ldots ,z_n)$. 
\end{proof}

\begin{proof}[{\bf Proof of Theorem~\ref{thm:sphere}}] Recall (\ref{eq:sphericalinvariance}) which asserts that $\X$ is invariant in distribution under the action of automorphisms of $\S^2$. By Lemma~\ref{lem:invdist}, it suffices to show that the density of points in $\X$ is of the form given in (\ref{eq:specialformofdensity}). We use the following well known matrix decomposition.

\noindent{\bf Schur decomposition:} Any diagonalizable matrix $M\in GL(n,\C)$
can be written as 
\begin{equation}\label{eq:schurdecomp}
  M = U(Z+T)U^*,
\end{equation}
where $U$ is unitary, $T$ is strictly upper triangular and $Z$ is
diagonal. Moreover the decomposition is almost unique, in the
following sense:

$M=V(W+S)V^*$ in addition to (\ref{eq:schurdecomp}), with $V,S,W$ being respectively  unitary, strictly upper triangular, and diagonal, if and only if
 the entries of $W$ are a permutation of the elements of $Z$, and if
 this permutation is identity, then 
 $V=U\Theta$ and $\Theta S \Theta^* =T$ for some $\Theta$ that 
 is both diagonal and unitary, that is, for 
 $\Theta$ of the form Diagonal$(e^{i\theta_1},\ldots ,e^{i\theta_n})$ .

Corresponding to this matrix decomposition (\ref{eq:schurdecomp}),
 Ginibre~\cite{gin} proved the following measure decomposition.

\noindent{\bf Ginibre's measure decomposition:} If $M$ is decomposed as in
(\ref{eq:schurdecomp}), with the elements of 
$Z$ in a uniformly randomly chosen order, then 
\begin{equation}
  \label{eq:ginibredecomp}
  \prod_{i,j} dm(M_{ij}) = \l(\prod_{i<j} |z_i-z_j|^2 \prod_k
  dm(z_k) \r)\l(\prod_{i<j} dm(T_{ij}) \r) d\nu(U)
\end{equation}
%where $\nu$ is a finite measure 
%given in (\ref{eq:haarunitarymoddiagunit}) 
where $\nu$ is the Haar measure on the on the unitary group $\U(n)$.
% such that $d\nu(U\Theta)=d\nu(U)$ for every
%diagonal unitary $\Theta$. 

Conditional on $G_1$, the matrix $M:=G_1^{-1}G_2$ has the density
\begin{equation*}
e^{-\Tr(M^*G_1^*G_1M)}|\det(G_1)|^{2n}
\end{equation*} 
with respect to the Lebesgue measure on $GL(n,\C)\subset \C^{n^2}$. 
From the measure decomposition  (\ref{eq:ginibredecomp}) we get the
density of $Z$, $T$, $U$, $G_1$ to be  
\begin{equation*} 
\l(\prod_{i<j} |z_i-z_j|^2 \prod_{k=1}^n
dm(z_k)\r) e^{-\Tr(G_1^*G_1(I+MM^*))} \Mid \det(G_1) 
\Mid^{2n} 
\end{equation*}
with respect to the measure  $d \nu(U) \prod_{i<j} dm(T_{ij}) \prod_{i,j}dm(G_1(i,j))$ (we have omitted constants entirely-they can be recovered at the end).
Thus the density of $Z$ is obtained by integrating over $T,U,G_1$. Now write $z_k=r_ke^{i\theta_k}$ so that  
$Z=\Theta R$ where $\Theta=\mb{Diagonal}(e^{i\theta_1},\ldots ,e^{i\theta_n})$ and $R=\mb{Diagonal}(r_1,\ldots ,r_n)$. Then
\be
 MM^*= U\Theta(R+\Theta^*T)(R+\Theta^*T)^*\Theta^*U^*. 
\ee 
As $\nu$ is the Haar measure, $d \nu(U\Theta)=d \nu(U)$. The elements of
$\Theta^*T$ are the same as elements of $T$, but multiplied by complex
numbers of absolute value $1$. Hence, $\Theta^*T$ has the same
``distribution'' as $T$. Thus replacing $U$ by $\Theta^*U$ and $T$ by
$\Theta^*T$ we see that the density of $Z$  
is of the form $V(r_1,\ldots ,r_n)\prodd_{i<j} |z_i-z_j|^2$. This is the
form of the density required to apply Lemma~\ref{lem:invdist}. Thus we
conclude that the eigenvalue density is 
\begin{equation}\label{eq:densityofX}
 \mb{Const.}\prodd_{i<j} |z_i-z_j|^2  \prod_{k=1}^n \frac{1}{
 (1+|z_k|^2)^{n+1}}.
\end{equation} 

To compute the constant, note that 
\begin{equation*}
\l\{ \sqrt{\frac{n}{\pi} \binom{n-1}{k}}\frac{z^k}{(1+|z|^2)^{\frac{n+1}{2}}} \r\}_{0\le k \le n-1}
\end{equation*}
is an orthonormal set. Projection on the Hilbert space generated by
 these functions  gives a determinantal process whose kernel
is as given in  the statement of the theorem. 
%Writing out the density shows
%that this is the same as the eigenvalue density that we have
% determined. Hence the constants must match.
\end{proof}

\section{Proof of Theorem~\ref{thm:general}}\label{sec:general}
We first find the coefficients in the power series expansion of $f_N$ prior to taking limits using the following lemma. Randomness plays no role here.
\begin{lemma}\lb{lem:blashke} Let $V$ be an $N\times N$ matrix  and define $f(z) = \frac{\det(zI+V)}{\det(I+zV^*)}$. Then
\be
 f^{(k)}(0) = \det(V) \summ_{\pi \in \Sym_k} \sgn(\pi) \prodd_{c\in \pi} \l[\Tr(V^{-|c|}) - \Tr(V^{*|c|}) \r]
\ee
where we write $c\in \pi$ to mean that $c$ is a cycle of $\pi$.
\end{lemma}
\begin{proof} Let $\chi(z)=\det(I+zV^{-1})$ and let $\psi(z)=\det(I+zV^*)$. Then $f(z)= \frac{\det(V)\chi(z)}{\psi(z)}$. Hence, 
\ben \label{eq:diffof1}
 f^{(k)}(z) = \det(V)\summ_{p=0}^k {k \choose p} \chi^{(k-p)}(z) \l(\frac{1}{\psi} \r)^{(p)}(z).
\een
First let us find the derivatives of $\psi$ and $\chi$. Let $V_{[j_1,\ldots ,j_k]}$ be the $k\times k$ matrix got from $V$ by deleting all rows and columns except the $j_1,\ldots ,j_k^{\mb{th}}$ ones.
\bes
 \psi^{(k)}(0) &=& k! \summ_{j_1<\ldots <j_k} \det\l(V_{[j_1,\ldots ,j_k]}^* \r) \\
&=& \summ_{(j_1,\ldots j_k)} \det\l(V_{[j_1,\ldots ,j_k]}^* \r) \hsp{5mm}  (\mb{ summand vanishes if }j_1=j_2) \\
 &=& \summ_{(j_1,\ldots j_k)} \summ_{\pi \in \Sym_k} \sgn(\pi) \prodd_{i=1}^k V_{j_i,j_{\pi(i)}}^* \\
 &=& \summ_{\pi \in \Sym_k} \sgn(\pi) \summ_{(j_1,\ldots j_k)} \prodd_{i=1}^k V_{j_i,j_{\pi(i)}}^*.
\ees
The inner sum factors over cycles of $\pi$. Let us write $c\in \pi$ to mean that $c$ is a cycle of $\pi$ and let $|c|$ denote the size of the cycle $c$.  Then we may write
\ben \lb{eq:diffpsi}
 \psi^{(k)}(0) = \summ_{\pi \in \Sym_k} \sgn(\pi) \prodd_{c \in \pi}\Tr(V^{*|c|}).
\een
Analogously, we have
\ben \lb{eq:diffchi}
 \chi^{(k)}(0) = \summ_{\pi \in \Sym_k} \sgn(\pi) \prodd_{c \in \pi}\Tr(V^{-|c|}).
\een
To compute the derivatives of $f$ using (\ref{eq:diffof1}), we need the derivatives of $\phi:=1/\psi$ at $0$. These will be given by the sequence that we shall provisionally call $\{b_k\}$. Set $b_0=1$ and for $k\ge 1$ define
\be
 b_k = \summ_{\pi \in \Sym_k} \sgn(\pi) \prodd_{c \in \pi}[-\Tr(V^{*|c|})].
\ee
Then for any $k\ge 1$ we calculate using (\ref{eq:diffpsi})
\bes
& & \summ_{j=0}^k {k \choose j} b_{k-j}\psi^{(j)}(0) \\
&=& \summ_{T\subset [k]} \l(\summ_{\pi\in \Sym(T^c) } \sgn(\pi) \prodd_{c\in \pi} \Tr(V^{*|c|}) \r) \l(\summ_{\pi\in \Sym(T) } \sgn(\pi) \prodd_{c\in \pi} [-\Tr(V^{*|c|})]  \r).
\ees
Fix a subset $T\subset [k]$. A permutation of $T$ and a permutation of $T^c$ together give a permutation of $[k]$. Let $\pi=c_1\ldots c_l$ be a permutation of $[k]$. Then it can arise from summands in which $T$ is a union (possibly empty)  of some of the cycles $c_i$s. Thus for $k\ge 1$
\bes
\summ_{j=0}^k {k \choose j} b_{k-j}\psi^{(j)}(0) &=&  \summ_{\stackrel{\pi \in \Sym_k}{\pi=c_1\ldots c_l}} \sgn(\pi) \summ_{B\subset [l]} \prodd_{i\in B} \Tr(V^{*|c_i|}) \prodd_{i\in B^c} \l[ -\Tr(V^{*|c_i|})\r] \\
&=& \summ_{\pi \in \Sym_k} \sgn(\pi) \prodd_{c\in \pi} \l[\Tr(V^{*|c|}) - \Tr(V^{*|c|})\r]\\
&=& 0.
\ees
However, the equation $\phi \cdot \psi =1$ implies that
\be
 \summ_{j=0}^k {k \choose j} \phi^{(k-j)}(0) \psi^{(j)}(0) = \l\{ \begin{array}{ll}
 1 & \mb{ if }k=0. \\ 0 & \mb{ if }k\not=0. \end{array} \r.
\ee
It is also clear that from these equations one may inductively recover $\phi^{(k)}(0)$ in terms of the derivatives of $\psi$. This shows that $\phi^{(k)}(0)=b_k$. That is
\ben \lb{eq:diffphi}
 \phi^{(k)}(0) = \summ_{\pi \in \Sym_k} \sgn(\pi) \prodd_{c\in \pi} \l[-\Tr \l(V^{*|c|} \r) \r].
\een
Now we return to the derivatives of $f$. From (\ref{eq:diffof1}), (\ref{eq:diffchi}) and (\ref{eq:diffphi}) we deduce that
\bes
& &f^{(k)}(0) = \det(V)\summ_{p=0}^k {k \choose p} \chi^{(k-p)}(0) \phi^{(p)}(0) \\ 
&=& \det(V) \summ_{T\subset [k]} \chi^{|T^c|}(0) \phi^{|T|}(0) \\  
&=&\det(V)\summ_{T\subset [k]} \l(\summ_{\pi\in \Sym(T^c) } \sgn(\pi) \prodd_{c\in \pi} \Tr(V^{-|c|}) \r) \l(\summ_{\pi\in \Sym(T) } \sgn(\pi) \prodd_{c\in \pi} [-\Tr(V^{*|c|})]  \r).
\ees
Just as before, a permutation of $T$ and a permutation of $T^c$ together give a permutation of $[k]$ and a permutation $\pi\in\Sym_k$ can arise from summands in which $T$ is a union (possibly empty) of some of the cycles of $\pi$. Therefore
\bes
f^{(k)}(0) &=& \det(V) \summ_{\stackrel{\pi \in \Sym_k}{\pi=c_1\ldots c_l}} \sgn(\pi) \summ_{B\subset [l]} \prodd_{i\in B} \Tr(V^{-|c_i|}) \prodd_{i\in B^c} \l[ -\Tr(V^{*|c_i|})\r] \\
&=& \det(V) \summ_{\pi \in \Sym_k} \sgn(\pi) \prodd_{c\in \pi} \l[\Tr(V^{-|c|}) - \Tr(V^{*|c|})\r].
\ees
\end{proof}

The probabilistic part of the theorem comes from the following lemma on Haar-distributed unitary matrices.
\begin{lemma}\label{lem:unitaries}
Let $U$ be an $N\times N$ random unitary matrix sampled from the Haar measure. Fix $n\ge 1$. After multiplication by $\sqrt{N}$, the first principal $n\times n$ sub-matrices of $U^p$, $p\ge 1$, converge in distribution to independent matrices with i.i.d. standard complex Gaussian entries. In symbols,
\be
 \sqrt{N}\l([U]_{i,j\le n}, [U^2]_{i,j\le n},\ldots  \r) \convd (G_1,G_2,\ldots )
\ee
where $G_i$ are independent $n\times n$ matrices with i.i.d. standard complex Gaussian entries. More precisely, any {\em finite} number of random variables $\sqrt{N}[U^p]_{i,j}$, $p\ge 1$, $i,j\le n$, converge in distribution to independent standard complex Gaussians.
\end{lemma}
In the literature, there are many results which are similar in spirit to Lemma~\ref{lem:unitaries}. For instance, Diaconis and Shahshahani~\cite{diasha94} computed (a slight mistake in that paper was corrected in Diaconis and Evans~\cite{de01}) showed that if $U$ is sampled from Haar measure on $\U(N)$, then $(\Tr(U), \Tr(U^2),\ldots ) \convd (g_1,\sqrt{2}g_2,\ldots)$.   Jiang~\cite{jiang}, answering a question of Diaconis, proved that  if $p_N,q_N$ are negligible compared to $\sqrt{N}$, then  the entries of the principal $p_N\times q_N$  submatrix of a unitary random matrix $U$ sampled from Haar measure on $\U(N)$, are approximately independent complex Gaussians. Our requirement is somewhere between the two. We need only submatrices of fixed size, but of all powers of $U$. We give a complete proof of Lemma~\ref{lem:unitaries} in section~\ref{sec:unitaries}.

\begin{proof}[{\bf Proof of Theorem~\ref{thm:general}}] 
Define $\f$ as in the statement of the theorem. Lemma~\ref{lem:blashke} asserts that 
\ben \lb{eq:anamika1}
 \f^{(k)}(0) = \det(V) \summ_{\pi \in \Sym_n} \sgn(\pi) \prodd_{c\in \pi} \l[\Tr(V^{-|c|}) - \Tr(V^{*|c|}) \r].
\een
We want to find the limit distribution of $\{\f^{(k)}(0):0\le k <\infty\}$. For this first let us consider $\Tr(V^{-p}) - \Tr(V^{*p})$ for $p\ge 1$. Setting $P=[P_1 : P_2]$ and $Q^*=[Q_1^*: Q_2^*]$ where $P_1,Q_1^*$ are $N\times n$ matrices, from (\ref{eq:lawofV}) we get
\ben \lb{eq:twomatrices}
 V^{-1} = P_1A^{-1}Q_1 + P_2Q_2 \mb{ and } V^* = P_1A^*Q_1 + P_2Q_2. 
\een
Then write
\bes
& & \Tr(V^{-p}) - \Tr(V^{*p}) \\
&=& \summ_{i_1,\ldots ,i_p} (V^{-1})_{i_1,i_2}\ldots(V^{-1})_{i_p,i_1} -   (V^*)_{i_1,i_2}\ldots(V^*)_{i_p,i_1} \\
&=& \summ_{i_1,\ldots ,i_p} \prodd_{j=1}^p (P_1A^{-1}Q_1 + P_2Q_2)_{i_j,i_{j+1}} - 
 \prodd_{j=1}^p (P_1A^*Q_1 + P_2Q_2)_{i_j,i_{j+1}}.
\ees
Here it is implied that $i_{p+1}=i_1$. Expand each of the products to get a sum of $2^p$ terms. Each of these terms is identified uniquely by an integer $0\le r\le p$ and a vector $\q=(q_1,\ldots ,q_r)$ of integers $q_1<q_2 <\ldots <q_r$ which are the values of $j$ for which we choose $(P_1A^{-1}Q_1 )_{i_j,i_{j+1}}$ (or $(P_1A^*Q_1)_{i_j,i_{j+1}}$), while for other $j$ we choose $(P_2Q_2)_{i_j,i_{j+1}}$ in both products. 

A most important observation is that all summands with $r=0$ cancel. What remains is
\ben \label{eq:idenidu} 
 \summ_{i_1,\ldots ,i_p} \summ_{r\ge 1} \summ_{\q}  \prodd_{j\notin \q} (P_2Q_2)_{i_j,i_{j+1}} \l(\prodd_{l=1}^r (P_1A^{-1}Q_1)_{i_{q_l},i_{q_l+1}} - \prodd_{l=1}^r (P_1A^*Q_1)_{i_{q_l},i_{q_l+1}}\r).
\een
We are using $\q$ to denote the vector $(q_1,\ldots ,q_r)$ as well as the set $\{q_1,\ldots ,q_r\}$ but this should not lead to any confusion. Now write for each $l=1,\ldots ,r$
\bes
 (P_1A^{-1}Q_1)_{i_{q_l},i_{q_l+1}} &=& \summ_{\alp_l=1}^n\summ_{\bet_l=1}^n (P_1)_{i_{q_l},\alp_l}(A^{-1})_{\alp_l,\bet_l}(Q_1)_{\bet_l,i_{q_l+1}}.\\
(P_1A^*Q_1)_{i_{q_l},i_{q_l+1}} &=& \summ_{\alp_l=1}^n\summ_{\bet_l=1}^n (P_1)_{i_{q_l},\alp_l}(A^*)_{\alp_l,\bet_l}(Q_1)_{\bet_l,i_{q_l+1}}.
\ees
Fix a choice of $r\ge 1$, $\q$ and $\alp_l,\bet_l$, $1\le l \le r$. Sum over $i_1,\ldots ,i_p$ in (\ref{eq:idenidu}). When we sum over $i_j$ for $q_1< j \le q_2$, in both the summands corresponding to $A^{-1}$ and $A^*$, we get a factor of (we have displayed only those factors that depend on $i_j$ for $q_1< j \le q_2$)
\be
\summ_{i_{q_1+1}\ldots i_{q_2}} (Q_1)_{\bet_1,i_{q_1+1}}\l[ \prodd_{j=q_1+1}^{q_2-1}(P_2Q_2)_{i_j,i_{j+1}}\r](P_1)_{i_{q_2},\alp_2}  = (Q_1(P_2Q_2)^{q_2-q_1-1}P_1)_{\bet_1,\alp_2}.
\ee
Similarly we sum over $i_j$ for $j$ between $q_l+1$ to $q_{l+1}$ for every $l$ (where $r+1=1$). Write $\lam_l=q_{l+1}-q_l-1$ for $l \le r-1$ and $\lam_r=r-q_r+q_1-1$. Then for a fixed value of $r\ge 1$, $\q$ and $\alp_l,\bet_l$, $1\le l \le r$, as we sum over all $i_j$s in (\ref{eq:idenidu}) we get 

\be
 \l[\prodd_{l=1}^r (A^{-1})_{\alp_l,\bet_l}-\prodd_{l=1}^r (A^{*})_{\alp_l,\bet_l}\r] \cdot \prodd_{l=1}^r (Q_1(P_2Q_2)^{\lam_l-1}P_1)_{\bet_l,\alp_{l+1}}.
\ee 
Any choice of $(\lam_1,\ldots ,\lam_r)$ comes from $p$ different choices of $\q$ (by cyclically rotating $(q_1,\ldots ,q_r)$). Therefore $\Tr(V^{-p}) - \Tr(V^{*p})$ is equal to
\ben \lb{eq:anamika2}
 p \summ_{\tiny{r\ge 1}} \summ_{\stackrel{\tiny{(\lam_1,\ldots ,\lam_r)}}{\tiny{\lam_i\ge 1, \sum \lam_i =p}}}  \summ_{\stackrel{\tiny{\alp_l,\bet_l}}{\tiny{l\le r}}}\l[\prodd_{l=1}^r (A^{-1})_{\alp_l,\bet_l}-\prodd_{l=1}^r (A^{*})_{\alp_l,\bet_l}\r] \cdot \prodd_{l=1}^r (Q_1(P_2Q_2)^{\lam_l-1}P_1)_{\bet_l,\alp_{l+1}}.
\een
As before, here $r+1=1$. 
Since $PQ$ has Haar distribution, from Lemma~\ref{lem:unitaries} and the assumption on $A$, we know that 
\be
 \sqrt{N}(A, Q_1(PQ)^0P_1, Q_1(PQ)^1P_1,Q_1(PQ)^2P_1,\ldots ) \convd (X_0,G_1,G_2,\ldots ).
\ee
where $G_i$, $i\ge 1$ are independent $n\times n$ matrices with i.i.d. standard complex Gaussian entries and independent of $X_0$. (Pre-multiplication by $Q_1$ and post-multiplication by $P_1$ serve to pick out the first $n\times n$ principal sub-matrix of $(QP)^m$). 

Now consider $Q_1(PQ)^mP_1$. For $m=0$, 
\be
 Q_1(PQ)^0P_1 =  Q_1(P_2Q_2)^0P_1
\ee
 Hence $\sqrt{N}Q_1(P_2Q_2)^0P_1 \convd G_1$. Next take $m=1$. Since  $PQ=P_1Q_1+P_2Q_2$, 
\be
Q_1(PQ)^1P_1 =  Q_1(P_2Q_2)^1P_1+ (Q_1P_1)^2.
\ee
From the $m=0$ case, we know that the second summand is $O_p(N^{-1})$, whence, $\sqrt{N}(Q_1(P_2Q_2)^0P_1,Q_1(P_2Q_2)^1P_1) \convd (G_1,G_2)$.
Continuing inductively, for any $m$, we get $Q_1(QP)^mP_1 = Q_1(P_2Q_2)^mP_1+ O_p(N^{-1})$. Thus
\be
 \sqrt{N}(A, [Q_1(P_2Q_2)^0P_1], [Q_1(P_2Q_2)^1P_1],\ldots ) \convd (X_0,G_1,G_2,\ldots )
\ee
in the sense that any finite subset of random variables on the left converge in distribution to the corresponding random variables on the right.

In equation (\ref{eq:anamika2}) divide each of the $r$ factors in the products inside the brackets by $\sqrt{N}$ and multiply each factor in the product outside the brackets by $\sqrt{N}$. $A^*$ itself converges to $0$ in probability and thus after dividing by $\sqrt{N}$, in the first product only $A^{-1}$ survives in the limit. Thus we get
\bes
& &\Tr(V^{-p}) - \Tr(V^{*p}) \\
&\convd& p \summ_{r\ge 1} \summ_{\stackrel{(\lam_1,\ldots ,\lam_r)}{\lam_i\ge 1, \sum \lam_i =p}}  \summ_{\alp_l,\bet_l,l\le r} \prodd_{l=1}^r X_0^{-1}(\alp_l,\bet_l) \cdot \prodd_{l=1}^r G_{\lam_l}(\bet_l,\alp_{l+1}) \\
&=& p \summ_{r\ge 1} \summ_{\stackrel{(\lam_1,\ldots ,\lam_r)}{\lam_i\ge 1, \sum \lam_i =p}} \summ_{\alp_l, l\le r}\prodd_{l=1}^r (X_0^{-1}G_{\lam_l})(\alp_l,\alp_{l+1}).
\ees
Use this in (\ref{eq:anamika1}) and observe that $\det(V)=\det(A)\det(P^*Q^*)=\det(A)e^{i\theta}$ where $\theta$ is uniform on $[0,2\pi]$ and independent of $A$. Absorb $e^{i\theta}$ into $X_0$ and denote $H_p=X_0^{-1}G_p$. Then we see that  $N^{\frac{n}{2}}\f^{(k)}(0)$ converges in distribution to (jointly for $k\ge 0$, of course) 
\ben \label{eq:prelimlimit}
 \det(X_0) \summ_{\pi \in \Sym_k} \sgn(\pi) \prodd_{c\in \pi} \l(|c|\summ_{r\ge 1} \summ_{\stackrel{\tiny{(\lam_1,\ldots ,\lam_r)}}{\tiny{\lam_i\ge 1, \sum \lam_i =|c|}}} \summ_{\alp_l,l\le r} \prodd_{l=1}^r (X_0^{-1}G_{\lam_l})(\alp_l,\alp_{l+1})\r)
\een
where $|c|$ is the number of elements in the cycle $c$. We must reduce this further. When we completely expand the products in (\ref{eq:prelimlimit}) we see that the right hand side is equal to (as usual $r_i+1=1$)
\ben \label{eq:prelim2limit}
\summ_{\stackrel{\pi \in \Sym_k}{\pi=c_1\ldots c_m}} \sgn(\pi)\prodd_{i=1}^m |c_i|  \summ_{\stackrel{r_i\ge 1}{1\le i\le m}} \summ_{\stackrel{(\lam_1^i,\ldots ,\lam_{r_i}^i)}{\lam_j^i\ge 1, \summ_j \lam_j^i =|c_i|}} \summ_{\stackrel{\alp_l^i\le n}{l\le r_i, i\le m}} \prodd_{i=1}^m \prodd_{j=1}^{r_i} H_{\lam_j^i}(\alp_{j}^i,\alp_{j+1}^i).
\een
The point is that many of the terms $\prodd_i \prodd_j H_{\lam_j^i}(\alp_{j}^i,\alp_{j+1}^i)$ can arise from more than one permutation $\pi$ and thus there is a lot of cancellation. This we investigate now.

Consider any term $\prodd_{l=1}^L H_{\mu_l}(s_l,t_l)$ where $L\ge 1$, $\mu_l\ge 1$ for each $l\le L$ and $1\le s_l,t_l \le n$.  We compute the coefficient of such a term in (\ref{eq:prelim2limit}). 

To organize the combinatorics that will emerge, for the term  $\prodd_{l=1}^L H_{\mu_l}(s_l,t_l)$ let us associate a directed multi-graph with edge-weights as follows. We assume that $\summ_{l=1}^L \mu_l=k$ as only such terms can appear in (\ref{eq:prelim2limit}). 

The graph will have vertices $\{1,2,\ldots ,n\}$. For each $l\le L$, put a directed edge  from $s_l$ to $t_l$  and give it weight $\mu_l$. Let us also put self loops with edge-weight $0$ at each vertex $v\in \{1,\ldots ,n\}\backslash \{s_l,t_l:l\le L\}$. Let us call this graph $\GG$ (depends on $\mu_l$, $s_l$ etc, of course, but it would be horrifying to include that dependence in the notation!). We group terms together by the graph they generate and find the total contribution for each graph.

The graph $\GG$ can arise from a term in (\ref{eq:prelim2limit}) only if the edges of $\GG$ can be partitioned into edge-disjoint directed cycles. Note that $\GG$ is a multi-graph and hence if $i\tends j$ occurs twice in $\GG$, the two instances will occur in two distinct cycles, but the cycles will be deemed disjoint. Also, a cycle may visit the same vertex more than once. 

 Furthermore, each such decomposition of $\GG$ into disjoint cycles corresponds to some (usually more than one) choice of the  permutation $\pi$ in (\ref{eq:prelim2limit}). Once $\pi$ is fixed, the numbers $r_i$ are just the sizes of cycles in this cycle decomposition of $\GG$ and $\lam_j^i$, $\alp_j^i$ are also determined. An example is given below to elucidate the matter.
\begin{example} Suppose $n=6$ and let $k=11$. Suppose we look at the term 
\ben \label{eq:term}
H_1(1,3)H_2(3,2)H_2(2,1)H_1(2,1)H_3(1,4)H_2(4,2).
\een
This term can actually arise in (\ref{eq:prelim2limit}) because the sum of the $\mu_l$s is equal to $k$ and the associated graph may be decomposed into disjoint cycles in two distinct ways: Firstly, as $\{(1,3,2), (1,4,2)\}$  and secondly, as $(1,3,2,1,4,2)$. 

The first case, $\{(1,3,2), (1,4,2)\}$,  can arise from any permutation $\pi \in \Sym_{11}$ that has two cycles of lengths $5$ and $6$ (these numbers come from adding the edge weights in each cycle). There are ${11 \choose 5}4!5!$ such permutations and they all have sign $-1$. Taking into account the weight $\prod |c_i|$ in (\ref{eq:prelim2limit}) , the contribution to the term (\ref{eq:term}) from all such permutations is $-11!$.

The second case, $(1,3,2,1,4,2)$, can arise from any $\pi\in \Sym_{11}$ that is itself a cycle of length $11$. There are $10!$ such permutations and they all have sign $+1$. Their total contribution is $+11!$. 
When put together, we see that the coefficient of (\ref{eq:term}) is zero. This is no coincidence and prepares the reader for what is stated next in general.
\end{example}

\noindent{\bf Claim:} If $\GG$ can be decomposed into disjoint cycles in more than one way, then the coefficient of the corresponding term in (\ref{eq:prelim2limit}) is zero. On the other hand, if $\GG$ can be decomposed in a unique way into $\ell$ disjoint cycles, then the coefficient of the corresponding term is $(-1)^{k-\ell} k!$. 

\noindent{\bf Proof of the Claim:} Suppose $\GG$ can be decomposed into disjoint cycles in more than one way. Then some vertex, say $1$ without losing generality, belongs to more than one cycle of $\GG$. Then in $\GG$, there are in-edges $i_1,\ldots ,i_M$ leading to $1$  and out-edges $j_1,\ldots ,j_M$ leading away from $1$, for some $M\ge 2$. In any decomposition of $\GG$ into cycles, we have the obvious matching of in-edges with out-edges, by associating to each in-edge the out-edge that follows it in the cycle. Consider any cycle decomposition of $\GG$, and suppose that in-edges $i_1,i_2$ are matched with $j_1,j_2$ respectively. We pair this cycle decomposition with a new cycle decomposition got by  switching the matches $i_1\tends j_1,i_2\tends j_2$ to $i_1 \tends j_2, i_2\tends j_1$  and leaving everything else intact. This leads to a pairing of all cycle decompositions of $\GG$. We show that the total contribution from each pair is zero.

One can go from one cycle decomposition to its pair by splitting a cycle into two cycles or merging two cycles into one. Let us take the first one among them to have cycle sizes $\theta_1,\ldots ,\theta_{\ell}$ and the second one to have cycle sizes $\theta_1+\theta_2,\theta_3,\ldots ,\theta_{\ell}$. Let the sums of edge-weights along cycles in the first decomposition be $w_1,w_2, \ldots ,w_{\ell}$ so that in the second cycle decomposition the sums of edge-weights of cycles are $w_1+w_2, w_3, \ldots ,w_{\ell}$. 

The permutations in $\Sym_k$ that respect the first cycle decomposition of $\GG$  are precisely those with $\ell$ cycles of sizes $w_1,\ldots ,w_{\ell}$. The number of such permutations is 
\ben \label{eq:somegarbage}
 \frac{k!}{\prodd_{i=1}^{\ell}w_i! } \prodd_{i=1}^{\ell}(w_i-1)!.
\een
Each of these comes with the weight $\prodd_{i=1}^{\ell} w_i$ in (\ref{eq:prelim2limit}), whence the total contribution of these terms is $(-1)^{k-\ell}k!$. 

The number of permutations that respect the second cycle decomposition of $\GG$ is\be
 \frac{k!}{(w_1+w_2)!\prodd_{i=3}^{\ell}w_i! } (w_1+w_2-1)!\prodd_{i=3}^{\ell}(w_i-1)!.
\ee
Each of these comes with the weight $(w_1+w_2)\prodd_{i=3}^{\ell} w_i$ in (\ref{eq:prelim2limit}), whence the total contribution of these terms is $(-1)^{k-\ell+1}k!$. 

Thus the two cycle decompositions exactly cancel each other out and it is seen that the total contribution is zero. This proves the first part of the claim.
 
For the second part, there is only one cycle decomposition by assumption, and the same calculations that led to (\ref{eq:somegarbage}) show that the coefficient is $(-1)^{k-\ell}k!$. This completes the proof of the claim.

\vspace{2mm}

Now consider a $\GG$ that has a unique cycle decomposition. Then that cycle decomposition may be regarded as a permutation $\tau \in \Sym_n$, where all those vertices that do not occur among $\{s_l,t_l\}$ are fixed points of $\tau$ (this is why we added self-loops to all these vertices when defining $\GG$). Recall that the edge-weights of these self-loops is $0$. It will be convenient to set $H_0=I_n$. Observe that $\sgn(\tau)=(-1)^{n-\ell}$. Then using the claim above to simplify (\ref{eq:prelim2limit}) we finally have
\ben \label{eq:limitfinally}
N^{\frac{n}{2}}\f^{(k)}(0) \convd (-1)^{k-n}k!\det(X_0) \summ_{\tau \in \Sym_n}\sgn(\tau) \summ_{\stackrel{(w_1,\ldots ,w_n)}{w_i\ge 0, \summ_i w_i =k}} \prodd_{i=1}^n (H_{w_i})_{i,\tau_i}.
\een
Here $H_p=X_0^{-1}G_p$ for $p\ge 1$ and $H_0=I_n$. If we forget the $(-1)^{k-n}$ factor, the right hand side of (\ref{eq:limitfinally}) is precisely $k!$ times the coefficient of $z^k$ in the power series expansion of $\det(X_0+zG_1+z^2G_2+\ldots )$ as may be seen by expanding the determinant as
\be
 \det(X_0)\summ_{\tau \in \Sym_n}\sgn(\tau) \prodd_{i=1}^n (I+zH_1+z^2H_2+\ldots )_{i,\tau(i)}.
\ee
The factor $(-1)^{k-n}$ is rendered irrelevant by changing $z$ to $-z$ and multiplying the whole function by $(-1)^n$. 

Thus we have proved that any the power series coefficients of $N^{\frac{n}{2}}\f_N$ converge jointly in distribution to the coefficients in the power series of $\det(G_0+zG_1+z^2G_2+\ldots )$, in the sense that any finite number of coefficients in the former converge jointly in distribution to the corresponding coefficients in the latter. This completes the proof of the theorem.
\end{proof}

\section{Hyperbolic ensembles}\label{sec:disk}
In this section we prove Theorem~\ref{thm:disk}. We shall make use of the following result of {\.Z}yczkowski and Sommers~\cite{zycsom}.
\begin{result}[{\.Z}yczkowski and Sommers(2000)]\label{thm:zycsom}
Let $U$ be an $(N+n)\times (N+n)$ unitary matrix sampled from Haar measure on $\U(N+n)$. Let $V$ be the $N\times N$ principal sub-matrix got by deleting the first $n$ rows and columns of $U$. Then the eigenvalues of $V$ form a  determinantal process in the unit disk $\D$ with kernel
\be
\Kdet_N^{(n)}(z,w) = \summ_{k=0}^{N-1} {-n-1\choose k}(-1)^k z^k \bar{w}^k
\ee
with respect to the reference measure $d\mu_n(z)=\frac{n}{\pi}(1-|z|^2)^{n-1}dm(z)$. 
\end{result}

\begin{remark} If $V$ is a matrix, then $\frac{\det(zI-V)}{\det(I-zV^*)}$ is just the Blaschke product of the eigenvalues $\{\lam_k\}$ of $V$. That is
\ben \label{eq:writingasblaschke}
 \frac{\det(zI-V)}{\det(I-zV^*)} = \prodd_{j=1}^k \frac{z-\lam_k}{1-z\bar{\lam}_k}.
\een
That being the case, the function $\f_N$ in Theorem~\ref{thm:disk} depends only on the eigenvalues and not the matrix $V$ that we choose. Why then, do we use the truncated unitary matrix of Result~\ref{thm:zycsom} instead of directly using the diagonal matrix whose entries are a determinantal process with the truncated kernels? It may indeed be possible to prove Theorem~\ref{thm:disk}, directly from the properties of determinantal processes without having to use Result~\ref{thm:zycsom}. However, that would involve proving a bevy of central limit theorems (of non-linear statistics) for determinantal processes that can substitute Lemma~\ref{lem:unitaries}. We do not know if that is easy. 

 The advantages of the truncated unitary matrix over the diagonal matrix of its eigenvalues are: (1) The former is invariant under left and right multiplication by unitary matrices (which allows us to apply the rather easy Lemma~\ref{lem:unitaries}). (2) The truncated unitary matrix has only $n$ nondeterministic singular values, even as the matrix size goes to infinity.  The cost is that we use the (far from trivial) result of {\.Z}yczkowski and Sommers but this has the positive value of forging a direct link between random matrices and random analytic functions.
\end{remark}

Applying Theorem~\ref{thm:general} to truncated unitary matrices we almost get Theorem~\ref{thm:disk}, but there is one snag. Theorem~\ref{thm:general} gives convergence of  coefficients in the power series whereas to deduce convergence of zeros, we need uniform convergence on compact sets. For instance, in the sequence $f_n(z)=n^nz^n$, all the power series coefficients converge to zero but $f_n(z)$ does not converge for any $z\not=0$. The following lemma, deduced directly from properties of determinantal point processes, will establish the required tightness, \'{a} priori.

\begin{lemma}\label{lem:tightness} Fix $n\ge 1$. Let $\{\lam_k:1\le k \le N\}$, be determinantal on the unit disk with kernel $\Kdet_N$ with respect to the measure $\mu_n$ as in Result~\ref{thm:zycsom}. Set
 \be
\f_N(z)=\prodd_{k=1}^N \frac{z-\lam_k}{1-z\bar{\lam}_k}.
\ee
Then for any compact subset $K\in \D$, the set $\{\f_N(z):z\in K\}$ is tight, uniformly in $N$.
\end{lemma}
We postpone the proof of Lemma~\ref{lem:tightness} to section~\ref{sec:tightness} and proceed to prove  Theorem~\ref{thm:disk} assuming the lemma.

\begin{proof}[{\bf Proof of Theorem~\ref{thm:disk}}] Let $U$ be an $(N+n)\times (N+n)$ unitary matrix  and write
\be
 U =  \l[\begin{array}{cc}
         A & C^* \\ B & V \end{array} \r]
\ee
where $A$ has size $n\times n$. By the unitarity of $U$, we have the following equations.
\be
 A^*A+B^*B=I_n \mb{ and } BB^*+VV^* = I_{N}.
\ee
As $BB^*$ and $B^*B$ have the same nonzero eigenvalues, it follows that $VV^*$ has the same eigenvalues as $A^*A$, except that it has $N-n$ more eigenvalues all equal to $1$. Thus there must exist unitary matrices $P,Q \in \U(N)$ such that
\ben \lb{eq:structureofV}
 V = P \l[\begin{array}{cc}
         A & 0 \\ 0 & I_{N-n} \end{array} \r] Q.
\een
Now suppose $U$ was sampled according to Haar measure on $\U(N+n)$. Then for any unitary matrices $P_0,Q_0 \in \U(N)$, we have
\be
\l[\begin{array}{cc}
         I_n & 0 \\ 0 & P_0 \end{array} \r]  \l[\begin{array}{cc}
         A & C^* \\ B & V \end{array} \r] \l[\begin{array}{cc}
         I_n & 0 \\ 0 & Q_0 \end{array} \r] \eqd  \l[\begin{array}{cc}
         A & C^* \\ B & V \end{array} \r]
\ee 
because Haar measure is invariant under left and right multiplication by group elements. This shows that $P_0VQ_0\eqd V$, which, together with (\ref{eq:structureofV}) implies that
\be
 V = P \l[\begin{array}{cc}
         A & 0 \\ 0 & I_{N-n} \end{array} \r] Q
\ee
where $P,Q,A$ are independent, $P,Q$ are distributed according to Haar measure on $\U(N)$ and $A$ is the principal $n\times n$ submatrix of an $(N+n)\times (N+n)$ unitary matrix.

Lemma~\ref{lem:unitaries} shows that $\sqrt{N}A \convd G_0$ where $G_0$ is an $n\times n$ matrix of i.i.d. standard complex Gaussians. Thus Theorem~\ref{thm:general} applies and we get
\be
N^{\frac{n}{2}} \frac{\det(zI-V)}{\det(I-zV^*)} \convd \det\l(\summ_{k=0}^{\infty} G_k z^k \r)
\ee
where all $G_k$, $k\ge 0$ are i.i.d. matrices of i.i.d. standard complex Gaussians.

This convergence is only in the sense of pointwise convergence of coefficients in the power series. But in case of truncated unitary matrices, Result~\ref{thm:zycsom} and Lemma~\ref{lem:tightness} together strengthen it to uniform convergence on compact subsets of $\D$. Therefore the zeros of $\f_N$ converge in distribution to the zeros of the limiting analytic function.

The upshot is that the point process of eigenvalues of $V$ (which are exactly the zeros of $\f_N$) converge in distribution to the zeros of $\det(G_0+zG_1+z^2G_2+\ldots)$. Use the result of {\.Z}yczkowski and Sommers and let $N\tends \infty$. The kernels $\Kdet_N$ increase (in the sense of operators, i.e., the associated Hilbert spaces increase) to the kernel $\Kdet(z,w)=(1-z\bar{w})^{-n-1}$. From the facts stated after definition~\ref{def:det} in the appendix,it follows that the determinantal process with kernel $\Kdet_N$ converges to the determinantal process with kernel $\Kdet$ and the proof is complete.
\end{proof}

\section{Proof of Lemma~\ref{lem:unitaries}}\label{sec:unitaries}
%Beginning of a supposedly simpler proof but not really any
% much simpler except perhaps psychologically. Omit for now.
%------------------------------------------------------------
\begin{comment}
It is well-known result that any submatrix (of fixed size) of a  Haar($\U(N)$) distributed unitary matrix, upon scaling by $\sqrt{N}$, converges in distribution to a matrix of i.i.d. standard complex Gaussians. See \cite{jiang} and the references therein. We shall use the following moments-version of this result. 
\begin{result}\label{res:unitaries} Let $U=((u_{i,j}))_{i,j\le N}$ be chosen from Haar measure on $\U(N)$. Fix $k\ge 1$ and consider  $1\le i(\ell),j(\ell),i'(\ell),j'(\ell) \le N$ for $1\le \ell \le k$. Then 
\be
N^k\E \l[\prodd_{\ell=1}^k u_{i(\ell),j(\ell)}  \prodd_{\ell=1}^{k'} \bar{u}_{i'(\ell),j'(\ell)} \r] 
\tends  \del_{k,k'} \summ_{\pi \in \Sym_k} \prodd_{\ell=1}^k \I{i(\ell)=i'(\pi_{\ell})}\I{j(\ell)=j'(\pi_{\ell})}. 
\ee
\end{result}
The matrix $M$ defined by $M(p,q)= \I{i(\ell)=i'(\pi_{\ell})}\I{j(\ell)=j'(\pi_{\ell})}$, is the covariance matrix of $\{g_{i(\ell),j(\ell)}:\ell \le k\}$ with $\{g_{i'(\ell),j'(\ell)}:\ell \le k'\}$, where $g_{i,j}$ are i.i.d. standard complex Gaussians. The right hand side in the statement of Result~\ref{res:unitaries} is precisely equal to $\mb{Per}(M)$. By the Wick formula for moments of complex Gaussians, 
\end{comment}
%-------------------------------------------------------------
% End of commented text

In proving Lemma~\ref{lem:unitaries}, we shall make use of  the following result on the joint moments of entries of a unitary matrix from the book of Nica and Speicher~\cite{nicspe}, page 381 (we state a weaker form suited to our purpose). 

\begin{result}\label{res:unitaries} Let $U=((u_{i,j}))_{i,j\le N}$ be chosen from Haar measure on $\U(N)$. Let $k\le N$ and fix $i(\ell),j(\ell),i'(\ell),j'(\ell)$ for $1\le \ell \le k$. Then 
\be
\E \l[\prodd_{\ell=1}^k u_{i(\ell),j(\ell)}  \prodd_{\ell=1}^k \bar{u}_{i'(\ell),j'(\ell)} \r] \\
= \summ_{\pi,\sig \in \Sym_k} \mb{Wg}(N,\pi \sig^{-1})\prodd_{\ell=1}^k \I{i(\ell)=i'(\pi_{\ell})}\I{j(\ell)=j'(\sig_{\ell})} 
\ee
where $\mb{Wg}$ (called ``Weingarten function'') has the property that as $N\tends \infty$, 
\be
 \mb{Wg}(N,\tau)  = \l\{ \begin{array}{ll}
                N^{-k} + O(N^{-k-1}) & \mb{ if }\tau=e \mb{ (``identity'' )}. \\
                O(N^{-k-1}) & \mb{ if }\tau \not= e. \end{array} \r.
\ee
%\frac{\phi(\tau)}{N^{2k-\mb{cyc}(\tau)}} + O\l(\frac{1}{N^{2k+2- \mb{cyc}(\tau)}}\r) \hsp{4mm} \mb{ as }N\tends \infty,
\end{result}

\begin{proof}[{\bf Proof of Lemma~\ref{lem:unitaries}}] We want to show that $\sqrt{N}(U^k)_{\alp,\bet}$, $k\ge 1$, $1\le \alp,\bet \le n$ converge (jointly) in distribution to independent standard complex Gaussians. To use the method of moments consider two finite products of these random variables
\be
 S = \prodd_{i=1}^m [(U^{k_i})_{\alp_i,\bet_i}]^{p_i} \hsp{1mm} \mb{ and }\hsp{1mm} T = \prodd_{i=1}^{m'} [(U^{k'_i})_{\alp'_i,\bet'_i}]^{p'_i}.
\ee
where $m,m',p_i,p'_i,k_i,k'_i \ge 1$ and $1\le \alp_i,\bet_i,\alp'_i,\bet'_i \le n$ are fixed. We want to find $\E[S\bar{T}]$ asymptotically as $N\tends \infty$.

The idea is simple-minded. We expand each $(U^k)_{\alp,\bet}$ as a sum of products of entries of $U$. Then we get a huge sum of products and we evaluate the expectation of each product using Result~\ref{res:unitaries}. Among the summands that do not vanish, most have the same contribution and the rest are negligible. We now delve into the details.

Let $\PP_k(\alp,\bet)$ denote all ``paths'' $\gam$ of length $k$ connecting $\alp$ to $\bet$. This just means that $\gam \in [N]^{k+1}$, $\gam(1)=\alp$ and $\gam(k+1)=\bet$. Then we write
\be
 (U^k)_{\alp,\bet} = \summ_{\gam \in \PP_k(\alp,\bet)} \prodd_{j=1}^{k} u_{\gam(j),\gam(j+1)}.
\ee
Expanding each factor in the definition of $S$ like this, we get
\be
 S = \summ_{\stackrel{\gam_i^{\ell}\in \PP_{k_i}(\alp_i,\bet_i)}{i\le m; \ell\le p_i}} \prodd_{i=1}^m \prodd_{\ell=1}^{p_i} \prodd_{j=1}^{k_i} u_{\gam_i^{\ell}(j),\gam_i^{\ell}(j+1)}.
\ee
In words, we are summing over a packet of $p_1$ paths of length $k_1$ from $\alp_1$ to $\bet_1$, a packet of $p_2$ paths of length $k_2$ from $\alp_2$ to $\bet_2$, etc. $T$ may similarly be expanded as
\be
 T = \summ_{\stackrel{\Gam_i^{\ell}\in \PP_{k_i'}(\alp_i',\bet_i')}{i\le m'; \ell\le p_i'}} \prodd_{i=1}^{m'} \prodd_{\ell=1}^{p_i'} \prodd_{j=1}^{k_i'} u_{\Gam_i^{\ell}(j),\Gam_i^{\ell}(j+1)}.
\ee
To evaluate $\E[S\bar{T}]$, for each pair of collections $\gam=\{\gam_i^{\ell}\}$ and $\Gam = \{\Gam_i^{\ell}\}$, we must find
\ben \label{eq:gamandGam}
 \E\l[ \l(\prodd_{i=1}^m \prodd_{\ell=1}^{p_i} \prodd_{j=1}^{k_i} u_{\gam_i^{\ell}(j),\gam_i^{\ell}(j+1)} \r) \l( \prodd_{i=1}^{m'} \prodd_{\ell=1}^{p_i'} \prodd_{j=1}^{k_i'} \bar{u}_{\Gam_i^{\ell}(j),\Gam_i^{\ell}(j+1)} \r)\r].
\een
Fix a collection of packets $\gam_i^{\ell}\in \PP_{k_i}(\alp_i,\bet_i)$. For which collections $\Gam_i^{\ell} \in \PP_{k_i'}(\alp_i',\bet_i')$  does (\ref{eq:gamandGam}) give a nonzero answer? For that to happen, the number of $u_{i,j}$s and the number of $\bar{u}_{i,j}$s inside the expectation must be the same (because $e^{i\theta}U\eqd U$ for any $\theta\in \R$). Assume that this is the case.
%This means that we must have $\summ_{i=1}^m p_ik_i = \summ_{i=1}^{m'} p_i'k_i' = r$ (say). Assume this happens. 

It will be convenient to write $\gam(i,\ell,j)$ in pace of $\gam_i^{\ell}(j)$. From Result~\ref{res:unitaries}, to get a nonzero answer in (\ref{eq:gamandGam}) we must have bijections 
\bes
\{(i,\ell,j):i\le m,\ell \le p_i ,1\le j\le k_i\} \stackrel{\pi}{\tends}\{(i,\ell,j):i\le m',\ell \le p_i' ,1\le j\le k_i'\}, \\
\{(i,\ell,j):i\le m,\ell \le p_i ,2\le j\le k_i+1\} \stackrel{\sig}{\tends}\{(i,\ell,j):i\le m',\ell \le p_i' ,2\le j\le k_i'+1\},
\ees
such that 
\bes
\l( \gam(i,\ell,j) \r)_{i\le m,\ell \le p_i ,1\le j\le k_i}  &=& \l( \Gam(\pi(i,\ell,j)) \r)_{i\le m,\ell \le p_i, 1\le j\le k_i}. \\
\l( \gam(i,\ell,j) \r)_{i\le m,\ell \le p_i ,2\le j\le k_i+1}  &=& \l( \Gam(\sig(i,\ell,j)) \r)_{i\le m,\ell \le p_i, 2\le j\le k_i+1}.
\ees 
And for each such pair of bijections $\pi,\sig$, we get a contribution of $\mb{Wg}(N,\pi\sig^{-1})$. 

Let us call the collection of packets $\gam$ {\bf typical}, if all the paths $\gam_i^{\ell}$ are pairwise disjoint (except possibly at the initial and final points) and also non self-intersecting (again, if $\alp_i=\bet_i$, the paths in packet $i$ intersect themselves, but only at the end points). 

If $\gam$ is typical, then it is clear that for $\Gam$ to yield a nonzero contribution, $\Gam$ must consist of exactly the same paths as $\gam$. This forces $k_i=k_i'$ and $p_i=p_i'$ and $\alp_i=\alp_i',\bet_i=\bet_i'$ for every $i$.  
%(Note: Although for ease of notation we index $\Gam$ and $\gam$ by $i,\ell$,  we are only regarding them as sets of paths. Collections of packets, which differ only in labelling are not counted multiple times).  
If this is so, then the only pairs of bijections $(\pi,\sig)$ that yield a non zero contribution are those for which
\begin{itemize}
\item $\pi = \sig$ (From the disjointness of the paths). 

\item $\pi$ permutes each packet of paths among itself. In particular there are $\prodd_{i=1}^k p_i!$ such permutations.
\end{itemize}
This shows that for a typical $\gam$, the expectation in (\ref{eq:gamandGam}) is equal to
\ben \label{eq:typical}
 \I{\Gam=\gam} \l( \prodd_{i=1}^m p_i! \r) \mb{Wg}(N,e).
\een
Here $\gam=\Gam$ means that the two sets of paths are the same.  
Now suppose $\gam$ is atypical. For any fixed $\gam$, typical or atypical, the number of $\Gam$ for which (\ref{eq:gamandGam}) is nonzero is clearly bounded uniformly by $m$ and $p_i,k_i$, $i\le m$. In particular it is independent of $N$. Therefore the expected value in (\ref{eq:gamandGam}) is bounded in absolute value by
\ben \label{eq:atypical}
 C \sup_{\tau}\mb{Wg}(N,\tau).
\een
Now for an atypical $\gam$, at least two of $\gam_i^{\ell}(j)$, $1\le i\le m$, $1\le \ell \le p_i$, $2\le j\le k_i$, must be equal (our definition of ``typical'' did not impose any condition on the initial and final points of the paths, which are anyway fixed throughout). Thus, if we set $r=p_1(k_1-1)+\ldots +p_m(k_m-1)$,  then it follows that the total number of atypical $\gam$ is less than $r^2N^{r-1}$. Since the total number of $\gam$ is precisely $N^r$, this also tells us that there are at least $N^r-r^2N^{r-1}$ typical $\gam$. Put these counts together with the contributions of each typical and atypical path, as given in (\ref{eq:typical}) and (\ref{eq:atypical}), respectively. Note that we get nonzero contribution from typical paths only if $S=T$. Also, the total number of factors in $S$ is $r+\sum p_i$ (this is the ``$k$'' in Result~\ref{res:unitaries}). Hence
\bes
\E[S\bar{T}] &=&  \I{S=T} N^r(1-O(1/N))\mb{Wg}(N,e) \prodd_{i=1}^m p_i!   + O(N^{r-1})  \sup_{\tau \in \Sym_{r+\sum p_i}}\mb{Wg}(N,\tau) \\
 &=& \I{S=T} N^{-\sum p_i} \l(\prodd_{i=1}^m p_i!\r) \l(1+O\l(\frac{1}{N}\r)\r)  
\ees
by virtue of the asymptotics of the Weingarten function, as given in Result~\ref{res:unitaries}.

The factor $N^{\sum p_i}$ is precisely compensated for, once we scale $(U^k)_{\alp,\bet}$ by $\sqrt{N}$, as in the statement of the lemma. Since the moments of standard complex Gaussian are easily seen to be $\E[g^p \bar{g}^q]=p!\I{p=q}$, we have shown that $\sqrt{N}(U^k)_{\alp,\bet}$, $k\ge 1$, $\alp,\bet \le n$,  converge to independent standard complex Gaussians.
\end{proof}

\section{Proof of Lemma~\ref{lem:tightness}}\label{sec:tightness}
We prove Lemma~\ref{lem:tightness} in this section. We shall make use of the following fact, which is a direct consequence of \cite{hkpv}, Theorem 26.

\begin{result}\label{res:absvalues} Fix $n>0$. Let $\{\lam_1,\ldots ,\lam_N\}$ be determinantal on the unit disk with kernel 
\be
\Kdet_N(z,w)=\summ_{k=0}^{N-1} {-n-1 \choose k}(-1)^k z^k \bar{w}^k
\ee
with respect to the background measure $d\mu_n(z)=\frac{n}{\pi}(1-|z|^2)^{n-1}dm(z)$. Then the set $\{|\lam_k|^2:1\le k\le N\}$ has the same distribution as $\{Y_k:0\le k\le N-1\}$, where $Y_k$s are independent random variables and $Y_k$ has distribution Beta($k+1,n$).
\end{result}

As a consequence of this result, it is very easy to see that $N^{\frac{n}{2}}\f_N(0)$ is tight. For, 
\bes
\E \l[|\f_N(0)|^2\r] &=& N^n \prodd_{k=0}^{N-1} \E[Y_k] \\
                     &=& N^n \prodd_{k=0}^{N-1} \frac{k+1}{n+k+1} \\
                     &=& N^n \frac{n!}{(N+1)(N+2)\ldots (N+n)} \\
                     &\tends& n!   
\ees
as $N\tends \infty$. For $z\not=0$ it is not as simple, because for finite $N$, the distribution of $\{\Mid \frac{z-\lam_k}{1-z\bar{\lam}_k} \Mid\}$ is not the same as that of a set of independent random variables. In the $N\tends \infty$ limit, it is, but is of no use to us.

\begin{proof}[{\bf Proof of Lemma~\ref{lem:tightness}}]  Write $\phi_z(\lam)=\frac{z-\lam}{1-z\bar{\lam}}$. Write
\be
 \summ_{k=1}^N \log |\phi_z(\lam_k)|^2 = -\summ_{k=1}^N (1-|\phi_z(\lam_k)|^2) + \summ_{k=1}^N h_z(\lam_k)
\ee
where $h_z(\lam) := \log|\phi_z(\lam)|^2+(1-|\phi_z(\lam)|^2)$. The lemma will be proved by showing that the following are tight (uniformly over $z$ in compact sets, as $N$ varies).
\begin{itemize}
\item $\summ_{k=1}^N h_z(\lam_k)$.

\item $-n \log N + \summ_{k=1}^N (1-|\phi_z(\lam_k)|^2)$.
\end{itemize}
Let us consider them one by one.

\begin{enumerate}
\item Consider $\summ_{k=1}^N h_z(\lam_k)$. For each $z\in \D$, the function $\phi_z$ maps the unit disk onto itself and the unit circle onto itself. Therefore, given a compact set $K$, we may find $T<1$ such that for any $z\in K$ and $|\lam|>T$ we have $|\phi_z(\lam)|^2>\frac{1}{2}$. From the power series expansion of $\log(1-x)$, we then get
\be
 |h_z(\lam)| < 2 (1-|\phi_z(\lam)|^2)^2 \hsp{3mm}\mb{ for }z\in K, |\lam|>T.
\ee
Observing that
\be
 1-|\phi_z(\lam)|^2 = \frac{(1-|z|^2)(1-|\lam|^2)}{|1-z\bar{\lam}|^2},
\ee
it follows that
\be
\given \summ_{k=1}^N h_z(\lam_k) \given \le 2\summ_{|\lam_k|<T} |h_z(\lam_k)|  + C \summ_{|\lam_k|>T} (1-|\lam_k|^2)^2
\ee
for a constant $C$ (does not depend on $N$ or $z$, as long as $z\in K$). The first summand is tight because the set $\{\lam_k: |\lam_k|<T\}$ converges to the set of points in $T\D$ in the limiting determinantal process. The second summand may be stochastically bounded by $\sum_{k=0}^{N-1} (1-Y_k)^2$ by Result~\ref{res:absvalues}. From the explicit distribution of $Y_k$s, we may compute the expected value of this sum as
\be
 \E\l[\summ_{k=0}^{N-1} (1-Y_k)^2 \r] = \summ_{k=0}^{N-1} \frac{n(n+1)}{(n+k+1)(n+k+2)}.
\ee
The random variables on the left are stochastically increasing in $N$ and hence, and hence, a uniform bound on the expectations shows tightness.
 
\item Consider $-n \log N + \summ_{k=1}^N (1-|\phi_z(\lam_k)|^2)$. We shall show tightness by proving that the expected value and variance of these random variables are bounded uniformly in $N$ (and $z\in K$). 

\noindent{\bf Expected value:} As always, it is simpler for us to deal with $|\lam_k|^2$. Therefore, write
\bes
1-|\phi_z(\lam)|^2 &=& \frac{(1-|z|^2)(1-|\lam|^2)}{|1-z\bar{\lam}|^2} \\
 &=& (1-|z|^2)(1-|\lam|^2)\summ_{p,q=0}^{\infty} z^p\bar{z}^q \lam^q\bar{\lam}^p.
\ees
Set $\lam=\lam_k$, sum over $k$ and take expectations. Any term with $p\not=q$ vanishes, because of rotation invariance of $\{\lam_k\}$. For $p=q$, we get terms with $|\lam_k|^2$ which may be replaced by independent Beta random variables by Result~\ref{res:absvalues}. Thus
\bes
& &\E\l[ \summ_{k=1}^N (1-|\phi_z(\lam_k)|^2) \r] \\
&=& (1-|z|^2)\summ_{p=0}^{\infty} |z|^{2p} \E\l[
\summ_{k=1}^N (1-|\lam_k|^2)|\lam_k|^{2p}  \r] \\
 &=&  (1-|z|^2)\summ_{p=0}^{\infty} |z|^{2p} \summ_{k=1}^N \frac{\mb{Beta}(k+p,n)-\mb{Beta}(k+p+1,n)}{\mb{Beta}(k,n)} \\
 &=& (1-|z|^2)\summ_{p=0}^{\infty}|z|^{2p} \summ_{k=1}^N \frac{n k(k+1)\ldots (k+n-1)}{(k+p)\ldots (k+p+n)} \\
 &=& (1-|z|^2)\summ_{p=0}^{\infty} |z|^{2p} \summ_{k=1}^N  \l( \frac{n}{k} + \mb{err}_k \r)
\ees
where $|\mb{err}_k| \le C(p) k^{-2}$ where $C(p)$ is at most a polynomial in $p$. Thus, it follows that
\bes
\E\l[ \summ_{k=1}^N (1-|\phi_z(\lam_k)|^2) \r] &=& (n\log N) (1-|z|^2) \l( \summ_{p=0}^{\infty}|z|^{2p} \r)+ O(1) \\
 &=& n\log N + O(1).
\ees
This is exactly what we wanted to show about the expected value.

\noindent{\bf Variance :} Now we want the variance of $X_N:=\summ_{k=1}^N (1-|\phi_z(\lam_k)|^2)$. No really new ideas are needed, only the calculations are more tedious. Expand $X_N$ in power series as before to get
\ben \label{eq:expandall}
E[X^2] =  \summ_{p,q,r,s\ge 0} z^{p+r}\bar{z}^{q+s}  \E\l[ (1-|\lam_k|^2)(1-|\lam_{\ell}|^2) \summ_{k,\ell=1}^N \lam_k^q \bar{\lam}_k^p \lam_{\ell}^s \bar{\lam}_{\ell}^r\r]. 
\een
All terms in which $p+r\not= q+s$ vanish, by rotation invariance. Fix $p,q,r,s$ so that $p+q=r+s$ and write the inner expectation as 
\be
 \E\l[\summ_{k=1}^N (1-|\lam_k|^2)^2|\lam_k|^{2p+2r} \r] + \E\l[\summ_{k\not= \ell} (1-|\lam_k|^2)(1-|\lam_{\ell}|^2) \lam_k^q \bar{\lam}_k^p \lam_{\ell}^s \bar{\lam}_{\ell}^r \r].
\ee
The first one is already a function of the absolute values of $\lam_k$s and hence replacing them by independent Beta random variables, we get (details are similar to those in computing the expectation) 
\bes
& &\E\l[\summ_{k=1}^N (1-|\lam_k|^2)^2|\lam_k|^{2p+2r} \r] =   \\
&=& \summ_{k=1}^N \frac{\mb{Beta}(k+p+r,n)-2\mb{Beta}(k+p+r+1,n)+\mb{Beta}(k+p+r+2,n)}
{\mb{Beta}(k,n)} \\
&=& n(n+1) \summ_{k=1}^N \frac{k(k+1)\ldots (k+n-1)}{(k+p+r) \ldots (k+p+r+n+1)}
\ees
which is bounded because the summand is of order $k^{-2}$. Of course, we shall have to sum over $p,q,r,s$, but it is clear that because of the factor of $|z|^{2p+2r}$ in (\ref{eq:expandall}), the total contribution to (\ref{eq:expandall}) from this summand (all terms with $k=\ell$) is bounded as $N\tends \infty$. 

It remains to consider the sum over $k\not=\ell$. The two point correlation is $\Kdet_N(\lam,\lam)\Kdet_N(\xi,\xi)-\Kdet_N(\lam,\xi)\Kdet_N(\xi,\lam)$. We consider
\bes
& &\E\l[\summ_{k\not= \ell} (1-|\lam_k|^2)(1-|\lam_{\ell}|^2) \lam_k^q \bar{\lam}_k^p \lam_{\ell}^s \bar{\lam}_{\ell}^r \r] \\
&=& \intt_{\D^2} (1-|\lam|^2)(1-|\xi|^2) \lam^q \bar{\lam}^p \xi^s \bar{\xi}^r \l(\Kdet_N(\lam,\lam)\Kdet_N(\xi,\xi)-|\Kdet_N(\lam,\xi)|^2\r).
\ees 
Consider the first summand, where we choose the term $\Kdet_N(\lam,\lam)\Kdet_N(\xi,\xi)$ inside the brackets. This survives only if $p=q$ and $r=s$ and it is easily seen that this term when summed over $p=q$ and $r=s$ in (\ref{eq:expandall}) will give exactly $\E[X]^2$. When we compute the variance of $X$, we shall subtract $\E[X]^2$ from $\E[X^2]$ and this term gets cancelled. 

Thus to show boundedness of the variance, we only need to show the boundedness of 
\ben \label{eq:doesthisneedaname}
\intt_{\D^2} (1-|\lam|^2)(1-|\xi|^2) \lam^q \bar{\lam}^p \xi^s \bar{\xi}^r \Kdet_N(\lam,\xi)\Kdet_N(\xi,\lam) d\mu_n(\lam)d\mu_n(\xi).
\een
Recall that 
\be
 \Kdet_N(\lam,\xi) = \summ_{j=0}^{N-1} C_j\lam^j \bar{\xi}^j.
\ee
where $C_j = {-n-1 \choose j}(-1)^j$. Also, for any $j$,
\be
 C_j \intt_{\D} |\lam|^{2j} d\mu_n(\lam) = 1.
\ee
Therefore the integral in (\ref{eq:doesthisneedaname}) is equal to
\bes
& & \summ_{i,j=0}^{N-1} C_iC_j\l(\intt_{\D}(1-|\lam|^2) \lam^{q+i} \bar{\lam}^{p+j} d\mu_n(\lam)\r)\l(\intt_{\D}(1-|\xi|^2) \xi^{s+j} \bar{\xi}^{r+i} d\mu_n(\xi) \r) \\
&=& \summ_{i,j=0}^{N-1} C_iC_j\del_{q+i,p+j}\del_{s+j,r+i} \l( \frac{1}{C_{p+j}}-\frac{1}{C_{p+j+1}}\r)\l( \frac{1}{C_{r+i}}-\frac{1}{C_{r+i+1}}\r).
\ees
We fixed $p,q,r,s$ such that $p-q=s-r$. Therefore, there are $N-|p-q|$ choices for $(i,j)$ which do not vanish (if $|p-q|\ge N$, there are no such terms). Without losing generality, let $p<q$, and write the above quantity as
\bes
& &\summ_{i=0}^{q-p-i} C_iC_{q-p-i} \l( \frac{1}{C_{q+i}}-\frac{1}{C_{q+i+1}}\r)\l( \frac{1}{C_{r+i}}-\frac{1}{C_{r+i+1}}\r) \\
&=& \summ_{i=0}^{N-|p-q|} \frac{n^2}{(i+1)(i+1+q-p)} + O(1)
\ees 
by writing out the expressions for $C_j$s. Thus, this term, when summed over $p,q,r,s$ yields a bounded quantity.

In summary, we wrote $E[X^2]$ as in (\ref{eq:expandall}). Terms with $k=\ell$ yielded a bounded quantity. Terms with $k\not=\ell$ were split into two sums. One of them is bounded while the other is exactly equal to $\E[X]^2$. Thus the variance is bounded as $N\tends \infty$.
\end{enumerate}
This completes the proof of Lemma~\ref{lem:tightness}.
\end{proof}

\section{Concluding remarks}
We record here two among several natural questions that arise from the considerations of this paper.
\begin{enumerate}
\item In Theorem~\ref{thm:disk}, the determinantal process exists for any positive real $n$, whereas the matrix analytic function makes sense only for integer values of $n$ (size of the matrix!). Is there a random zero set interpretation for hyperbolic determinantal processes for non-integer values of $n$?
 
\item Are there random matrices of the form given in Theorem~\ref{thm:general} for which we can calculate the exact distribution of eigenvalues? Recall that these are random matrices for which:  (1) The distribution is invariant under multiplication by unitary matrices. (2) The number of random singular values stays fixed even as the size of the matrix goes to infinity. 

If this can be done, then presumably we shall also get the distribution of zeros of $\det(X_0+zG_1+z^2G_2+\ldots)$ (we say "presumably" because we have not proved uniform convergence of $\f_N$ on compact sets except in the the special case of truncated unitary matrices). 

\end{enumerate}

\noindent{\bf Acknowledgements:} I thank B\'{a}lint Vir\'{a}g for asking me the question of finding determinantal processes on the sphere, Yuval Peres for innumerable illuminating discussions and for great encouragement throughout the project, and Mikhail Sodin for suggesting the use of Wick calculus which led to a verification of two-point correlations in \cite{krithesis} and convinced us that Theorem~\ref{thm:disk} must be true. I am greatly indebted to Brian Rider for pointing out the paper of {\.Z}yczkowski and Sommers which finally enabled me to prove Theorem~\ref{thm:disk}.

\section{Appendix: Determinantal point processes}\label{sec:appendix2} 
We give a brief introduction to determinantal processes, strictly limited to the context of this paper. More details, as well as proofs, may be found in the surveys~\cite{soshnikov00} or \cite{hkpv}. 

 Let $\Ome$ be a region in the plane and let $p$ be a positive continuous function on $\Ome$. Define the measure $\mu$ by  $d\mu(z) = p(z)dm(z)$. A simple point process on $\Ome$ is a random measure on $\Ome$ that takes values in counting measures on $\Ome$ and gives finite measure to compact sets. If the number of points in any compact set has exponential tails, then the distribution of the point process is determined by its correlation functions (joint intensities) with respect to $\mu$,
\be
\rho_k(z_1,\ldots ,z_k)  = \lim_{\eps\downarrow 0} \frac{\P[\X \mb{ has points in each of }D(z_i,\eps), 1\le i\le k]}{ \prodd_{i=1}^k \mu(D(z_i,\eps)) }
\ee
for any $k\ge 1$ and any $z_1,\ldots ,z_k\in \Ome$. There is also an integral version of this definition which is more appropriate in more general situations. Joint intensities need not exist in general.  

Consider the space $L^2(\Ome,\mu)$ and its subspace $\H$ consisting of holomorphic functions. It is a fact that  $\H$ is a closed subspace. For, suppose $f_n\in \H$ and $f_n\tends f$ in $L^2(\mu)$. Then for any $z\in \Ome$, consider a disk $D(z,r)$ contained entirely in $\Ome$. Convergence in $L^2$ and absolute continuity of $\mu$ shows that $f_n$s converge to $f$ in $L^2(m)$ on the annulus $D(z,r) \backslash D(z,\frac{r}{2})$. Apply Cauchy's integral formula to infer that $f_n \tends f$ uniformly on $D(z,\frac{r}{4})$. Therefore $f$ is itself analytic on $\Ome$. This shows that $\H$ is a Hilbert space. This reasoning also shows that for any $z\in \Ome$, the evaluation $f\tends f(z)$ is a bounded linear functional on $\H$.

Switch notations and let $\H$ denote any closed subspace of $L^2(\mu)$ consisting of holomorphic functions (not necessarily all holomorphic functions). Then the evaluation $f\tends f(z)$ is a bounded linear functional on $\H$. As a consequence, if  $\{\psi_n\}_{n\ge 1}$ is any orthonormal basis of $\H$, then the series
\be
\Kdet(z,w)=\sum \psi_n(z)\bar{\psi}_n(w),
\ee
does converge and $\Kdet$ (called the reproducing kernel of $\H$) is independent of the choice of the basis.
The integral operator on $L^2(\mu)$ defined by
\be
 \K f(z)  = \intt_{\Ome} f(w) K(z,w) d\mu(w)
\ee
is precisely the projection operator on $L^2(\mu)$ onto the subspace $\H$. 

\begin{definition}\label{def:det} Let $\H$ be a closed subspace of $L^2(\mu)$ consisting of analytic functions, and let $\Kdet$ be the reproducing kernel of $\H$. Then a point process $\X$ on $\Ome$ with joint intensities given by
\be
 \rho_k(z_1,\ldots ,z_k) = \det \l(\Kdet(z_i,z_j) \r)_{i,j\le k}
\ee
does exist and is called the determinantal point process on $\Ome$ with kernel $\Kdet$ with respect to the measure $\mu$. We also say that $\H$ is the associated Hilbert space.
\end{definition}
We just state a few facts regarding these processes.
\begin{enumerate}
\item The number of points in $\X$ is almost surely equal to the dimension of $\H$.
\item If subspaces $\H_n$ increase to $\H$, then the corresponding determinantal processes converge in distribution, i.e., $\X_n \convd \X$. This was used tacitly in deducing Theorem~\ref{thm:disk} from Result~\ref{thm:zycsom}.
\item If $\H=\mb{span}\{1,z,\ldots ,z^{n-1}\}$, then by writing out the density  (which is just $\frac{1}{n!}\rho_n(\cdot)$), one sees that the vector of $n$ points of the the process has density proportional to $\prodd_{i<j}|z_i-z_j|^2 \prodd_{i=1}^n p(z_i)$ with respect to the Lebesgue measure on $\Ome^n$.
\end{enumerate}

\bibliography{../Latex_commonfiles/mybibliography}
\end{document}